\documentclass[sn-mathphys]{sn-jnl}

\catcode`\|=12\relax

\usepackage{amssymb}
\usepackage{amsmath}
\usepackage{amsfonts}
\usepackage{upgreek}
\usepackage{siunitx}
\usepackage{listings}
\usepackage{multirow}
\usepackage{multicol}
\usepackage{bigdelim}
\usepackage{color}
\usepackage{stmaryrd}
\usepackage{soul}
\usepackage[symbol]{footmisc}
\usepackage{booktabs}
\usepackage[section]{placeins}
\usepackage[capitalize]{cleveref}
\usepackage{siunitx}

\usepackage{subcaption}
\usepackage{tablefootnote}

\usepackage[outdir=figures/]{epstopdf}
\usepackage{graphicx}
\usepackage{graphics}
\usepackage{varwidth}
\usepackage{adjustbox}
\usepackage{tikz}
\usepackage{pgfplots}
\usetikzlibrary{arrows,matrix,positioning,fit}
\usetikzlibrary{shapes,positioning}
\usetikzlibrary{backgrounds}
\pgfplotsset{compat=1.14}
\usepgfplotslibrary{colorbrewer}
\usepgfplotslibrary{patchplots}
\usepgfplotslibrary[colorbrewer]
\usetikzlibrary{pgfplots.colorbrewer}
\usetikzlibrary[pgfplots.colorbrewer]
\usepgfplotslibrary{fillbetween}
\usepgfplotslibrary{units}
\usetikzlibrary{spy}
\usepackage{pgfplotstable}
\graphicspath{{figures/}}
\usetikzlibrary{external}

\jyear{2021}

\theoremstyle{thmstyleone}%

\theoremstyle{thmstyletwo}%

\theoremstyle{thmstylethree}%

\newcommand\px[2]{\frac{\partial #1}{\partial {#2}}}

\usepackage{natbib}
\begin{document}

	\title[Nonlinear $ p $-multigrid preconditioner]{Nonlinear $ p $-multigrid preconditioner for implicit time integration of compressible Navier--Stokes equations}
	
	\author*[1]{\fnm{Lai} \sur{Wang}}\email{laiwang@tamu.edu}
	\author[2]{\fnm{Will} \sur{Trojak}\email{wtrojak@imperial.ac.uk}}
	\author[3]{\fnm{Freddie} \sur{Witherden}}\email{fdw@tamu.edu}
	\author[1]{\fnm{Antony} \sur{Jameson}}\email{antony.jameson@tamu.edu}
	\affil[1]{\orgdiv{Department of Aerospace Engineering}, \orgname{Texas A\&M University}}
	\affil[2]{ \orgdiv{Department of Aeronautics}, \orgname{Imperial College London}}
	\affil[3]{ \orgdiv{Department of Ocean Engineering}, \orgname{Texas A\&M University}}
	
\abstract{
     Within the framework of $ p $-adaptive flux reconstruction,
 we aim to construct efficient polynomial multigrid ($p$MG) preconditioners for implicit time integration of the Navier--Stokes equations using Jacobian-free Newton--Krylov (JFNK) methods.
 We hypothesise that in pseudo transient continuation (PTC), as the residual drops, the frequency of error modes that dictates the convergence rate gets higher and higher.
 We apply  nonlinear $p$MG solvers to stiff steady problems at low Mach number ($\mathrm{Ma}=10^{-3}$) to verify our hypothesis.
 It is demonstrated that once the residual drops by a few orders of magnitude, improved smoothing on intermediate $ p $-sublevels will not only maintain the stability of $ p $MG at large time steps but also improve the convergence rate. For the unsteady Navier--Stokes equations, we elaborate how to construct nonlinear preconditioners using pseudo transient continuation for the matrix-free generalized	minimal residual (GMRES) method used in explicit first stage, singly diagonally implicit Runge--Kutta (ESDIRK) methods, and linearly implicit Rosenbrock--Wanner (ROW) methods. Given that at each time step the initial guess in the nonlinear solver is not distant from the converged solution, we recommend a two-level $p\{p_0\text{-}p_0/2\} $ or even $ p\{p_0\text{-}(p_0-1)\} $ $p$-hierarchy for optimal efficiency with a matrix-based smoother on the coarser level based on our hypothesis. It is demonstrated that insufficient smoothing on intermediate $p$-sublevels will deteriorate the performance of $p$MG preconditioner greatly. The nonlinear $p$MG preconditioner in this framework is found to be effective in reducing computational cost, as well as reducing the dimension of Krylov subspace for stiff systems arising from high-aspect-ratio elements and low Mach numbers. Specifically, the JFNK-$ p $MG technique is demonstrated to be more than 5 times faster than $ p $MG nonlinear solvers for unsteady problems. Compared to the EJ preconditioner, the $ p $MG preconditioner can make ESDIRK and ROW methods up to 2 times faster for low-Mach-number flow and up to 1.5 times faster for highly anisotropic meshes. Moreover, the $p$MG preconditioner can reduce the dimension of Krylov subspace by one order of magnitude. With a $ p $MG preconditioner, ROW methods are consistently more efficient than ESDIRK methods.

}

    
\maketitle	
\section{Introduction}\label{intro}
    In order to accurately predict turbulent separated flows, scale-resolving simulations are needed since  traditional Reynolds-averaged Navier--Stokes (RANS) modeling approaches are not sufficently accurate. Scale resolving simulations with high-order methods have been demonstrated to be more efficient and accurate for these types of problems when compared to second-order methods that are widely used in industry~\cite{vermeire2017utility, jia2019evaluation}.
    In particular, implicit large eddy simulation (ILES)~\cite{boris1990large} with high-order methods have become  popular in the computational fluid dynamics (CFD) community over the years.
    Results from under-resolved simulations with high-order methods have been demonstrated to be reliable and accurate~\cite{Gassner2013}.
    The turn around time of LES of flows at moderate Reynolds numbers has been significantly reduced to days~\cite{vermeire2017utility, wang2017towards}.
    In spite of the recent progress, the Courant--Friedrichs--Lewy (CFL) condition has been a major restriction in tackling  high Reynolds number flows or flow systems that have strong source terms, such as chemical reaction flows. Unconditionally $ L $-stable implicit time integrators, however, enable performing time integration with time steps where the  CFL number is $ \gg 1$, thus significantly decreasing the number of time steps. With efficient low-storage nonlinear solvers as well as preconditioners, implicit methods is of great potential to accelerate scale resolving simulations at high Reynolds numbers and they are also more robust in  dealing with low-quality meshes.
    
     It has been an ongoing topic to seek efficient nonlinear solvers that achieve high convergence rates with minimised memory usage for massive parallel simulations of unsteady flows.
     Dating back to the last century, Jameson~\cite{jameson1991time} proposed dual time stepping, also widely known as pseudo transient continuation (PTC), which was coupled to a geometric multigrid ($h$MG) method~\cite{jameson1983solution} to solve time dependent problems.
     Since then, different solution strategies of $h$MG have been established~\cite{mavriplis1998multigrid, venkatakrishnan1995agglomeration, katz2009multicloud} in the literature to improve the efficiency of $h$MG methods. As high-order methods emerged, the polynomial multigrid method was proposed by R{\o}nquist and Patera~\cite{ronquist1987spectral} as an alternative, and is analogous to $h$MG in that a hierarchy of resolutions is built up using a hierarchy of polynomial degrees~\cite{bassi2003numerical,helenbrook2003analysis,fidkowski2005p,luo2006p,liang2009p}. 
     Alternatively, an inexact Newton--Krylov method, usually coupled with pseudo transient continuation, is another highly popular nonlinear solver used in the CFD community. In practice, for unsteady problems it is typical to employ the Jacobian--free Newton--Krylov (JFNK) method where the matrix-vector product in the construction of the Krylov subspace is approximated by a finite difference approximation~\cite{knoll2004jacobian,wang2020comparison}.
     To some extent,  in the past decade,  $p$MG or JFNK accelerated high-order methods have been successfully applied to LES~\cite{loppi2019locally, wang2020comparison, wang2020dynamically}. However, the solution strategy in~\cite{loppi2019locally} of using $p$MG only suffers from low convergence rates for stiff systems and the element-Jacobi preconditioner used in~\cite{wang2020comparison, wang2020dynamically} lacks robustness with limited dimension of Krylov subspace for large time steps. And even when a matrix-based approach is used, a large number for the Krylov subspace is recommended in~\cite{bassi2015linearly}.

     Mavriplis et al.~\cite{mavriplis1998multigrid} demonstrated that a Newton--Kyrlov method with multigrid precondtioners works better than a multigrid nonlinear solver for finite volume methods.
     More recently, Shahbazi et al.~\cite{shahbazi2009multigrid} showed that within the DG framework, $hp$MG preconditioned GMRES works better than $hp$MG methods alone. For unsteady simulations of stiff systems, we are hence motivated to pursue  efficient $p$MG preconditioned JFNK methods for implicit time integration.
     Regardless of whether $p$MG is employed as a nonlinear solver or preconditioner, one indispensable question is how to efficiently configure $p$MG for optimal performance. 
     Fidkowski et al.~\cite{fidkowski2005p} developed the element-line preconditioner for $p$MG with a two-level $\{p_0\text{-}(p_0-1)\}$ polynomial hierarchy in DG applied to the Navier--Stokes equations, where $p_0$ is the order of the finest $p$-sublevel. 
     This choice of polynomial hierarchy hides the impact of errors on $p$-resolution less than $p_0-1 $ since a direct solver is used on this lower $p$-sublevel. Persson and Peraire~\cite{persson2008newton} also focused on a two-level $p$MG preconditioner. They found that when the coarse $p$-sublevel employed $\mathbb{P}^1$ polynomials, the convergence of the $p$MG preconditioner could be almost Mach number independent when using an ILU0 smoother.
     Moreover,  Shahbazi et al.~\cite{shahbazi2009multigrid}  concluded that  $\mathbb{P}^0$ should be accurately solved for Euler equations and $ \mathbb{P}^1 $ should be sufficiently solved for Navier--Stokes equations, where sufficient smoothing is achieved via employing a saw-tooth like $ hp $-V-cycle starting from the coarsest level. 
     Even though in their work a dense polynomial hierarchy, $ \{p_0\text{-}(p_0-1)\text{-}\dots\text{-}1\} $, was used,   investigation of the intermediate $p$-sublevels was not made.
     Wang and Yu~\cite{wang2019p, wang2019implicitpmg} recommended to employ $\{p_0\text{-}p_0/2\text{-}p_0/4\}$ as the polynomial hierarchy for 3-level $p$MG nonlinear solvers based on preliminary numerical results, which worked well for the Euler equations but was not satisfying for both steady and unsteady Navier--Stokes equations. Note that $\mathbb{P}^0$ was employed at the coarsest $p$-sublevel and only the EJ smoother was used. Franciolini et al.~\cite{franciolini2020efficient} recently coupled $ p $MG and $p$-adaptation for scale-resolving turbulence simulation using discontinuous Galerkin (DG) methods. In their work, the strategy was to use an element-Jacobi preconditioned and matrix-free GMRES on the finest level and matrix-based GMRES with different preconditioners on coarser levels for a dense polynomial hierarchy, which  was observed to be faster and less memory-consuming than a matrix-based method with single $p$-resolution.
     
     In general, the effect of intermediate $p$-sublevels on convergence for either the $p$MG nonlinear solver or preconditioner has been widely overlooked in the literature, with analysis mainly performed for steady problems. We remark that for unsteady problems, the solution at time step $ t^n $ is not very distant from that at time step $ t^{n+1} $ even for implicit time integration, which is quite different from steady problems where the initial guess is typically significantly distinct from the converged solution. These overlooked factors motivates us to conduct the present study.  
    
    First, we analyse the performance of the $p$MG method as a nonlinear solver coupled to flux reconstruction~\cite{huynh2007flux,wang2009unifying,vincent2011new} for steady solutions of the Navier--Stokes equations. We investigate the effect of smoother strength  at $ p $-sublevels and how this can affect the rate of convergence.  Experiments have been done to demonstrate our hypothesis that in PTC, as the residual drops the frequency of the error modes that dictates the convergence rate increases. Therefore,  it is possible to use the PTC residual as an indicator to gradually switch smoother on intermediate $ p $-sublevel from element-Jacobi (EJ) to MBNK to accelerate convergence at the expense of small memory overhead. Second, we formulate a general framework of constructing nonlinear $p$MG preconditioners for Krylov subspace methods via pseudo transient continuation, focusing on GMRES in this work. Since the solution at time step $t^{n}$ is comparatively close to that at time step $ t^{n+1} $ for unsteady simulation, the dictating error modes are theorised to be on the higher intermediate $p$-sublevels. Hence, we investigate employing a strong smoother at intermediate $p$-sublevel (such as $ p_0/2 $ or $p_0-1$) to address these dominant error modes. 
    We have observed that the $p$MG preconditioner can significantly speed up the simulation for stiff systems due to the presence of large aspect ratio meshes and flows in the low-Mach-number regime. Moreover, the $ p $MG preconditioner can decrease the dimension of the Krylov subspace by one order of magnitude compared to the EJ preconditioner, which is vitally important for linearly implicit Rosenbrock--Wanner (ROW) methods~\cite{wanner1996solving}.  In addition, our analysis couples the $p$MG solver/preconditioner with the $p$-adaptive flux reconstruction method developed by Wang and Yu~\cite{wang2020dynamically}. 
    
    The remainder of this paper is organized as follows.~In \cref{sec:background}, we briefly review explicit first stage, singly diagonally implicit Runge--Kutta (ESDIRK) methods and Rosenbrock--Wanner methods. Using the pseudo transient continuation methodology, we elaborate the iterative methods that will be used for inexact Newton/Krylov-subspace methods in this study. In \cref{sec:pmgpreconditioner}, we first introduce the generic nonlinear $p$MG method and then formulate the framework of $p$MG preconditioner for GMRES method used in both ESDIRK and ROW. And we also address the numerical strategy of using $p$MG preconditioner for $p$-adaptive flux reconstruction method.  In \cref{sec:num_steady}, we present numerical experiments using $p$MG for steady problems and \cref{sec:num_unsteady} demonstrates the benefits of using $ p $MG preconditioners for stiff systems. Finally, the conclusions of this study are drawn in \cref{sec:conclusion} and we provide some insights for future research directions, specifically on applying $p$MG preconditioned JFNK for LES on GPU hardware.
\section{Background}\label{sec:background}
    \subsection{Implicit time integration}
    The compressible Navier--Stokes equations can be written in the general form of a conservation equation as:
    \begin{equation}\label{key}
        \px{\boldsymbol{q}}{t} + \nabla\cdot\boldsymbol{f} = 0.
    \end{equation}
    In this work we use the flux reconstruction (FR) method to discretise spatial derivatives. The FR method --- originally introduced by Huynh~\cite{huynh2007flux} --- has been extended to several element topologies as well as advection-diffusion equations, and we refer the reader to the works of Hynuh~\cite{huynh2007flux},  Wang and Gao~\cite{wang2009unifying}, Vincent et al.~\cite{vincent2011new}, Wang and Yu~\cite{wang2018compact}, and references therein for more details. 
    
    In this study two time integration methods are considered for unsteady calculations, namely the explicit-first-stage singly diagonally implicit Runge--Kutta (ESDIRK) methods and linearly implicit Rosenbrock--Wanner methods. Considering the integration from time step $n$ to $n+1$, the ESDIRK method reads as 
    \begin{equation}\label{ESDIRK2}
        \begin{cases}
            \boldsymbol{q}^{n+1} = \boldsymbol{q}^{n}+\Delta t\,
            \sum_{i=1}^{s}b_i\boldsymbol{R}(\boldsymbol{q}_i),\\
            \boldsymbol{q}_1 = \boldsymbol{q}^n,\\
            \boldsymbol{q}_i = \Delta t \, a_{ii}\, \boldsymbol{R}(\boldsymbol{q}_i)+\boldsymbol{q}^n+\Delta t\, \sum_{j=1}^{i-1}a_{ij}\boldsymbol{R}(\boldsymbol{q}_j),\,i=2,\dots,s,	
        \end{cases}	
    \end{equation}
    where $s$ is the number of stages and $\boldsymbol{R}(\boldsymbol{q})$ refers to the spatial discretization of $\nabla\cdot\boldsymbol{f}$. The ROW method can be considered  \emph{linearly implicit}, this refers to the linearization of traditional DIRK schemes in the ROW method such that for each stage one only needs to solve a linear system --- instead of a nonlinear one.  The general form of Rosenbrock methods applied to a system of conservative equations can be written as
    \begin{equation}\label{Rosenbrock_General}
        \begin{cases}
            \boldsymbol{q}^{n+1} = \boldsymbol{q}^n+\sum_{j=1}^{s}m_j\boldsymbol{Y}_j,\\	
            \left(\frac{\boldsymbol{I}}{a_{ii} \Delta t}-\frac{\partial \boldsymbol{R}}{\partial\boldsymbol{q}}\right)^n\boldsymbol{Y}_i = \boldsymbol{R}\left(\boldsymbol{q}^n+\sum_{j=1}^{i-1}a_{ij}\boldsymbol{Y}_j\right)+\frac{1}{\Delta t}\sum_{j=1}^{i-1}c_{ij}\boldsymbol{Y}_j,\ i=1,2,\ldots,s.
        \end{cases}
    \end{equation}
    Details of a series of ESDIRK and ROW methods from second order to fourth order can be found in~\cite{wang2020comparison}. When using ESDIRK, it is typical to use an inexact Newton's method to approximately solve the linearised system; for example, when using GMRES a reduction in the residual by one to two orders of magnitude is often sufficient. We should remark at this point that the stability and accuracy of ROW is solely contingent upon whether the linear solver can sufficiently reduce the residual, with several orders of magnitude typically being required.
    
    Second- and fourth- order ESDIRK and ROW methods, i.e., ESDIRK2~\cite{kennedy2016diagonally}, ESDIRK4~\cite{bijl2002implicit}, ROW2 (ROS2PR)~\cite{rang2014analysis}, and ROW4 (RODASP)~\cite{steinebach1995order}, are considered in this work.
    
\subsection{Iterative methods}
    To solve the nonlinear equations defined in \cref{ESDIRK2}, we can reformulate them as 
    \begin{equation} \label{eq:nonlinear}
        \boldsymbol{F}(\boldsymbol{q}_i) = 0, 
    \end{equation}
    where 
    \begin{equation}\label{SDIRK_newton}
        \boldsymbol{F}(\boldsymbol{q}_{i}) = \left(\frac{1}{a_{ii}\,\Delta t}\boldsymbol{q}_i -\boldsymbol{R}(\boldsymbol{q}_i)\right)-\frac{1}{a_{ii}\,\Delta t}\left(\boldsymbol{q}^n+\Delta t\sum_{j=1}^{i-1}a_{ij}\boldsymbol{R}(\boldsymbol{q}_j)\right).
    \end{equation}
    The PCT method, proposed by Jameson~\cite{jameson1991time}, is widely used to solve this system by finding a steady state solution in pseudo-time. This is achieved by introducing    a pseudo-time derivative into \cref{eq:nonlinear}, which can then be iteratively solved using a backward Euler method as
    \begin{equation}\label{pseudo_transient}
        \px{\boldsymbol{q}_i}{\tau} \approx\frac{\boldsymbol{q}^{k+1}_i-\boldsymbol{q}^{k}_i}{\Delta \tau}=-\boldsymbol{F}(\boldsymbol{q}^{k+1}_i).
    \end{equation}
    Herein, $k$ is the iteration step for the pseudo-transient continuation. \cref{pseudo_transient} can then be linearized as 
    \begin{equation}\label{linearized_pseudo_transient_final}
        \left(\frac{\boldsymbol{I}}{\Delta \tau} +\frac{\boldsymbol{I}}{a_{ii}\Delta t}-\frac{\partial \boldsymbol{R}}{\partial\boldsymbol{q}}\right)^k\Delta \boldsymbol{q}^{k}_i= -\boldsymbol{F}(\boldsymbol{q}^{k}_i).
    \end{equation}
    
    As the pseudo time is marched forward, the series of linear equations given by \cref{linearized_pseudo_transient_final} are successively solved until convergence, i.e. $\boldsymbol{q}^{k}_i\to\boldsymbol{q}_i^{n+1}$ as $k\to\infty$. 
    If matrix-free GMRES is employed as the linear solver, the matrix-vector product used in building the Krylov subspace is approximated as 
    \begin{equation}
        \left(\frac{\boldsymbol{I}}{\Delta \tau} +\frac{\boldsymbol{I}}{a_{ii} \Delta t}-\frac{\partial \boldsymbol{R}}{\partial\boldsymbol{q}}\right)\boldsymbol{X} = \boldsymbol{AX} \approx  \left(\frac{\boldsymbol{I}}{\Delta \tau}+\frac{\boldsymbol{I}}{a_{ii}\Delta t}\right)\boldsymbol{X} - \frac{\boldsymbol{R}(\boldsymbol{q}+\varepsilon \boldsymbol{X})-\boldsymbol{R}(\boldsymbol{q})}{\varepsilon},
    \end{equation}
    where a fixed value $\varepsilon = \num{1e-6}$ is used in this work.  We will use a fixed $\Delta\tau$, which if not otherwise stated will be set as $\Delta \tau/\Delta t = 1$ for unsteady problems. In the literature, when an explicit smoothing technique is used, a locally adpative pseudo time stepping, for example Loppi et al.~\cite{loppi2019locally}, can be more rewarding.
   
    It is well known that the performance of GMRES is strongly dependent on preconditioning. Many preconditioning methods are available, among which the element-Jacobi preconditioner is one of the simplest in the context of high-order spatial discretizations. For the element-Jacobi precondtioner, one uses the inverse of the diagonal blocks of the  matrix $\boldsymbol{A}$ to get the preconditioned vector $\boldsymbol{Y}$ as 
    \begin{equation}
        \boldsymbol{Y}=\mathrm{diag}(\boldsymbol{A})^{-1}\boldsymbol{X} = \boldsymbol{D^{-1}}\boldsymbol{X}.
    \end{equation}
    This assumes that $\mathrm{diag}(\boldsymbol{A})^{-1}$ is a good surrogate for $\boldsymbol{A}^{-1}$, which for preconditioning of the Navier--Stokes equations is acceptable when the stiffness is mild. Throughout this work the preconditioner matrix was evaluated at the beginning of each physical time step for JFNK methods and no Jacobian-freezing techniques were employed in this study.
   
    When considering steady problems, $\Delta t$ related terms in \cref{SDIRK_newton} are dropped, and we can use the successive evolution relaxation (SER) algorithm~\cite{mulder1985experiments} to set $\Delta\tau$ as
    \begin{equation}
        \Delta\tau^{0} = \Delta\tau_\mathrm{init},\quad
        \Delta \tau^{k+1} = \min\left(\Delta \tau^{k}\frac{\|\boldsymbol{F}^{k+1}\|_{L_2}}{\|\boldsymbol{F}^k\|_{L_2}}, \Delta \tau_\mathrm{max} \right).
    \end{equation}
    Where in the limit as $\Delta\tau \to \infty$, Newton's method is recovered. 
    
    Due to the $L$-stability of implicit time integration methods, the CFL number is now solely an indicator of the relative scaling of $\Delta t$ to the characteristic time scale of the simulation. There is still active discussion of a reliable and meaningful definition of CFL number for high-order methods on curved grids for the Navier--Stokes equations; we opt, therefore, to simply provide the minimum element size as a reference for interested readers, while our implementation runs with a dimensionless setup.
    
\section{$p$MG solver and $p$MG preconditioner}\label{sec:pmgpreconditioner}

\subsection{The nonlinear $p$MG solver}
   For completeness, we present the general procedure of a two-level $ p $MG method here. The hierarchy of the polynomial degrees is $ \{p_0\text{-}p_1\} $. Herein, subscripts `0' or `1' are the first and second $p$-levels.  This two-level hierarchy leads us to the following nonlinear system at each level, which can be expressed as 
    \begin{subequations}
        \begin{align}
            \boldsymbol{F}_{p_0}(\boldsymbol{q}_{p_0})-\boldsymbol{S}_{p_0} &= 0,\label{lvl_p_0}\\
            \boldsymbol{F}_{p_1}(\boldsymbol{q}_{p_1})-\boldsymbol{S}_{p_1} &= 0,\label{lvl_p_1}
        \end{align}
    \end{subequations}
    where $ {\boldsymbol{F}}$ is defined as \cref{SDIRK_newton} for unsteady problems and $\boldsymbol{F}=\boldsymbol{R}$ for steady problems. Here $ \boldsymbol{S}_{p_0} $ and $\boldsymbol{S}_{p_1}$ are referred to as forcing terms, noting that by definition
    \begin{equation}
        \boldsymbol{S}_{p_0} = 0.
    \end{equation} 
   
    \begin{figure}[tbhp]
   	    \centering   	
   	    \includegraphics[width=0.8\textwidth]{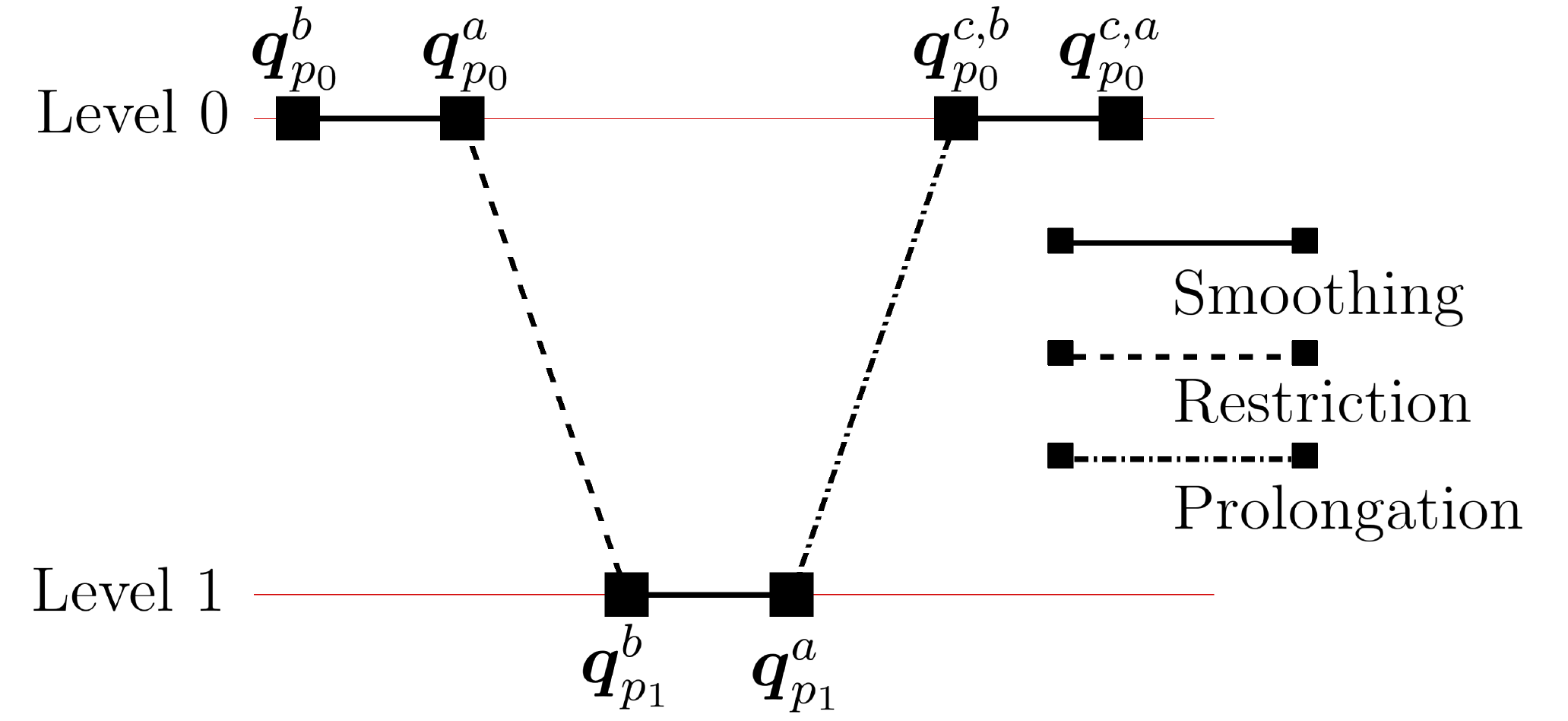}   	
   	    \caption{ Illustration of a 2-level V-cycle.}	
   	    \label{pmg-illustration}	
    \end{figure}
   
    The procedure for a typical $p$MG approach~\cite{fidkowski2005p,luo2006p,liang2009p} is then shown diagrammatically in \cref{pmg-illustration}, and has the following steps:
   
    \begin{description}
        \item[$p_0$ Smoothing.] Before smoothing, the initial value of $ \boldsymbol{q} $ at the first level is $ \boldsymbol{q}_{p_0}^b $. \cref{lvl_p_0} is then smoothed using an approach such as the element-Jacobi smoother for a number of steps. The primitive variables after smoothing is expressed as $ \boldsymbol{q}^{a}_{p_0}$. The defect at the first level 
   	    \begin{equation}
   	        \boldsymbol{d}_{p_0} = \boldsymbol{S}_{p_0}-\boldsymbol{F}_{p_0}(\boldsymbol{q}^{a}_{p_0}).
   	    \end{equation}
        \item[Restriction.] Restrict $\boldsymbol{q}^{a}_{p_0}$ at level $p_0$ to the lower $p_1$ level as 
   	    \begin{equation}
   	        \boldsymbol{q}_{p_1}^{b} = \mathbf{\Gamma}_{p_1}^{p_0}\boldsymbol{q}_{p_0}^{a},
   	    \end{equation}
        where $\mathbf{\Gamma}_{p_1}^{p_0}$ is a restriction operator from $p_0$ to $p_1$. The results, $\boldsymbol{q}_{p_1}^{b}$, is then the initial solution at the second level. Calculate the forcing term at the second level as 
   	    \begin{equation}
   	        \boldsymbol{S}_{p_1} = \boldsymbol{F}_{p_1}(\boldsymbol{q}_{p_1}^{b}) + \mathbf{\Gamma}_{p_1}^{p_0} \boldsymbol{d}_{p_0}
   	    \end{equation}
        \item[$p_1$ Smoothing.] Smooth \cref{lvl_p_1} to obtain the up-to-date solution $\boldsymbol{q}_{p_1}^a$. At this coarsest level we use either the element-Jacobi smoother or the matrix-based Newton--Krylov smoother.
   	
        \item[Prolongation.] Prolongate the change to the solution at level $p_1$ to produce a correction to the solution at the finest $p_0$ level. This can be written as
   	    \begin{equation}
        	\boldsymbol{q}_{p_0}^{c,b}= \boldsymbol{q}_{p_0}^{a}+\mathbf{\Pi}_{p_0}^{p_1}\boldsymbol{C}_{p_1}
   	    \end{equation}
   	    where $ \boldsymbol{C}_{p_1} = \boldsymbol{q}_{p_1}^{a}-\boldsymbol{q}_{p_1}^{b} $ and $\mathbf{\Pi}_{p_0}^{p_1}$ is the prolongation operator. 
   	    
        \item[$p_0$ Post-smoothing.] Post-smooth \cref{lvl_p_0} using the element-Jacobi smoother for a few steps with starting value $\boldsymbol{q}_{p_0}^{c,b}$ to obtain the smoothed solution $\boldsymbol{q}_{p_0}^{c,a}$ at the finest level. The result, $\boldsymbol{q}_{p_0}^{c,a}$, is the final solution, $\boldsymbol{q}_{p_0}$, after one V-cycle.
    \end{description}
    
    Throughout we use the notation that `b' means before smoothing, `a' means after smoothing, and `c' means corrected solution at the current $p$-level. In this work, the restriction operator, $\mathbf{\Gamma}_{p_1}^{p_0}$,  is  defined as the $L_2$ projection from the $p_0$ space to the coarser $p_1$ space, which reads~\cite{fidkowski2005p,luo2006p}
    \begin{equation}\label{projection}
        \mathbf{\Gamma}_{p_1}^{p_0} = (M^{p_1})^{-1}M^{p_1p_0}.
    \end{equation} 
    Herein, $M^p$ is the mass matrix for basis in a polynomial space $\mathbb{P}^p$, defined as $ M_{ij}^{p} = \int l_i^pl_j^p$. The second mass matrix, $M^{p_0p_1}$, is for a space $\mathbb{P}^{p_0}\times\mathbb{P}^{p_1}$ and is defined as $ M_{ij}^{p_1p_0} = \int l_i^{p_1}l_j^{p_0}$. Note that $ l_i^p $ is the $ i $-th Lagrange polynomial basis for the standard element which has a polynomial of degree $ p $. 
    The prolongation operator $ \mathbf{\Pi}_{p_0}^{p_1} $ is defined as
    \begin{equation}\label{prolongation}
        \mathbf{\Pi}_{p_0}^{p_1} = (M^{p_0})^{-1}M^{p_0p_1}.
    \end{equation}
    Since the Lagrange polynomials are not orthogonal, these operators are generally dense matrices.
    In terms of smoothing, we use notation $n\{n_0\text{-}n_1\}$ for a two-level $ p $MG which means we use $n_0$ iterations for both pre- and post- smoothing on the $p_0$-sublevel and $n_1$ iterations for $p_1$-sublevel.  
     
\subsection{Nonlinear $p$MG preconditioner}
    If the linear equation to be solved is expressed as
    \begin{equation}\label{linear-system}
        \boldsymbol{AX}=\boldsymbol{b},
    \end{equation}
    the following equation
    \begin{equation}\label{left-preconditioning}
        \boldsymbol{Y}=\boldsymbol{P}^{-1}\boldsymbol{X},
    \end{equation}
    describes a preconditioning procedure in general, where $ \boldsymbol{Y} $ is the preconditioned vector, $ \boldsymbol{X} $ is the vector to be preconditioned, and $ \boldsymbol{P} $ is the preconditioning matrix.   If $ \boldsymbol{P^{-1}=A^{-1} }$, one can get the ideal preconditioned vector $ \boldsymbol{Y} $ (the solution vector).
    Linear preconditioners has mainly focused on how to get a better estimation of $ \boldsymbol{A^{-1}}  $.
    However, in preconditioning, we essentially want to get a better estimation of $ \boldsymbol{A^{-1}X} $. This offers a possibility to introduce the pseudo transient continuation technique. 
     Given that $ \boldsymbol{A} $ can be written as 
    \begin{equation}
        \boldsymbol{A} = \left(\frac{\boldsymbol{I}}{\Delta \tau}+\frac{\boldsymbol{I}}{a_{ii}\,\Delta t}-\frac{\partial \boldsymbol{R}}{\partial \boldsymbol{q}}\right),
    \end{equation}
    for ESDIRK methods. A linear system from \cref{left-preconditioning} can be written as 
    \begin{equation}\label{linear-left-pre}
        \left(\frac{\boldsymbol{I}}{\Delta \tau}+\frac{\boldsymbol{I}}{a_{ii}\,\Delta t}-\frac{\partial \boldsymbol{R}}{\partial \boldsymbol{q}}\right)\boldsymbol{Y}=\boldsymbol{X}.
    \end{equation}
    It can be reorganized as 
    \begin{equation}\label{eq:ptc_linear_preconditioner}
        0=\boldsymbol{X}+\frac{\partial \boldsymbol{R}}{\partial \boldsymbol{q}}\boldsymbol{Y}-\left(\frac{\boldsymbol{I}}{\Delta \tau}+\frac{\boldsymbol{I}}{a_{ii}\,\Delta t}\right)\boldsymbol{Y}.
    \end{equation}
    As we know the approximation $\boldsymbol{R}(\boldsymbol{q}+\boldsymbol{Y}) - \boldsymbol{R}(\boldsymbol{q}) \approx \px{R}{q}\boldsymbol{Y} $, we further reorganize \cref{eq:ptc_linear_preconditioner} as 
    \begin{equation}\label{eq:ptc_nolinear_preconditioner}
        \begin{split}
            0=\boldsymbol{X}+\boldsymbol{R}(\boldsymbol{q}+\boldsymbol{Y})- \left(\frac{\boldsymbol{I}}{\Delta \tau}+\frac{\boldsymbol{I}}{a_{ii}\,\Delta t}\right)(\boldsymbol{q}+\boldsymbol{Y}) 
         - \boldsymbol{R}(\boldsymbol{q}) +\left(\frac{\boldsymbol{I}}{\Delta \tau}+\frac{\boldsymbol{I}}{a_{ii}\,\Delta t}\right)\boldsymbol{q}.
        \end{split}
    \end{equation}   
    By letting $ \hat{\boldsymbol{q}} = \boldsymbol{q}+\boldsymbol{Y}$, the following nonlinear equation is to be solved
    \begin{equation}\label{nonlinear-precon}
        \begin{split}
            0=\boldsymbol{X}+\boldsymbol{R}(\hat{\boldsymbol{q}})-\left(\frac{\boldsymbol{I}}{\Delta \tau}+\frac{\boldsymbol{I}}{a_{ii}\,\Delta t}\right)\hat{\boldsymbol{q}}-\boldsymbol{R}(\boldsymbol{q})+\left(\frac{\boldsymbol{I}}{\Delta \tau}+\frac{\boldsymbol{I}}{a_{ii}\,\Delta t}\right)\boldsymbol{q}
        \end{split}  
    \end{equation}
    in the pseudo transient continuation preconditioner. Introducing  another pseudo time stepping with step size $\Delta\tau^*$ using the backward Euler method yields
    \begin{equation}\label{eq:esdirk-ptc-marching-precon}
        \begin{split}
            \frac{\hat{\boldsymbol{q}}^{m+1}-\hat{\boldsymbol{q}}^{m}}{\Delta \tau^*}=\boldsymbol{X}+\boldsymbol{R}(\hat{\boldsymbol{q}}^{m+1})-\left(\frac{\boldsymbol{I}}{\alpha\Delta \tau}+\frac{\boldsymbol{I}}{a_{ii}\,\Delta t}\right)\hat{\boldsymbol{q}}^{m+1}-\boldsymbol{R}(\boldsymbol{q})+\left(\frac{\boldsymbol{I}}{\Delta \tau}+\frac{\boldsymbol{I}}{a_{ii}\,\Delta t}\right)\boldsymbol{q},
        \end{split} 
    \end{equation}
    where $m$ is used to differentiate the iteration number for the preconditioner from other iteration numbers. As we march in pseudo time $\tau^*$, $\boldsymbol{Y}=\hat{\boldsymbol{q} }-\boldsymbol{q}$ gives a better preconditioned vector $\boldsymbol{Y}$, and for all preconditioning techniques this is the primary objective. For instance, a more accurate approximation of the inverse of the Jacobian matrix will also end up with a closer approximation of $\boldsymbol{Y}=\boldsymbol{A}^{-1}\boldsymbol{X}$.  
    We remark that the pseudo transient continuation preconditioner has made no assumption about the governing equations and neither do any operation depends on the specific form of the governing equation. Hence it can be applied to any type of problems if pseudo transient continuation can be introduced. Even though we use ESDIRK as an example to formulate the general framework using pseudo transient continuation. One can also employ the same approach for ROW methods via dropping terms which contain $ \Delta \tau $ in \cref{eq:esdirk-ptc-marching-precon}. To avoid confusion, we always set $\Delta\tau^* = \Delta\tau$ if they are both present in the simulation.

    If we use an element-Jacobi method to iteratively solve the  Eq.~\eqref{eq:esdirk-ptc-marching-precon} after linearization, one pseudo time stepping makes this preconditioner recover  the element-Jacobi preconditioner. Moreover, we can directly employ the $ p $MG method described in  the previous subsection to solve this nonlinear system. We refer to this approach as the nonlinear $ p $MG preconditioner to distinguish from other linear ones. In the rest of this paper, we use $ p $MG preconditioner since there is no ambiguity. There are more choices of different smoothers than the element-Jacobi preconditioner such as the ILU0 smoother used in~\cite{persson2008newton}, the element-line preconditioner developed in~\cite{fidkowski2005p}, which are more effective than classic ones, such as the Gauss-Seidel method as well as EJ. However, we will limit our discussion to element-Jacobi and matrix-based Newton--Krylov (MBNK) smoothers. For the matrix-based Newton--Krylov smoother, we only use $ \mathrm{tol}_\mathrm{gmres}^\mathrm{r} = 10^{-1}$ as the relative tolerance for GMRES and set the maximum number of GMRES iterations as 5, which we find are sufficient enough for numerical experiments considered in this study. In the matrix-based Newton--Krylov smoother, we use ILU0 as the preconditioner for each block when domain decomposition is used for parallel computing.

\subsection{Coupling with $ p $-adaptation}\label{coupling_adaptation}
    
    We have demonstrated that $ p $-adaptation is very effective to save both computational time and total number of degrees of freedom for scale-resolving simulations with proper load balancing techniques~\cite{wang2020dynamically}. The feature-based $ p $-adaptation employs the spectral decay smoothness indicator to adjust the polynomial degrees in the flow field. The spectral decay smoothness indicator has been successful used to detect trouble cells for shock-capturing~\cite{persson2006sub}. It is defined as 
    \begin{equation}
    \eta_k = \frac{\left\|s_p-s_{p-1}\right\|_{L^2}}{\left\|s_p\right\|_{L^2}}
    \end{equation}
    for an element $ k $. The polynomial degree $ p $ of an element in the flow field is  $ p\in[p_\mathrm{min},p_\mathrm{max}] $, where $ p_\mathrm{min} $ and $ p_\mathrm{max} $ are the minimum and maximum polynomial degree, respectively. If not specifically mentioned, we choose the momentum in $ x $ direction, namely $ \rho u $, as the quantity for smoothness indicator calculation; $\nu_\mathrm{max} = 0.2  $ and $ \nu_\mathrm{min} = 0.001 $ are employed. We are going to increase the local polynomial degree by 1 when $ \eta_k >\nu_\mathrm{max}\eta_\mathrm{max} $ and decrease it by 1 when  $ \eta_k <\nu_\mathrm{min}\eta_\mathrm{max} $, where $ \eta_\mathrm{max} $ is the maximum smoothness indicator in the flow field.  We refer readers  to~\cite{wang2020dynamically} for detailed information of $ p $-adaptation that we have developed. In this study, when we refer to $ p $-adaptive $ \mathbb{P}^k $ FR, the maximum polynomial degree in the flow field is $ p_\mathrm{max} = k $ and the minimum polynomial degree is always $ p_\mathrm{min} =1 $.
    
    The necessity of a hierarchy of different $ p $-resolutions in $ p $MG methods does not stop us from employing them for implicit time integration with $ p$-adaptive flux reconstruction methods, in which numerical solutions are approximated by non-uniform polynomial degrees. Different from the hierarchy of $ p $-uniform $ p $MG method illustrated in Figure~\ref{pmg-illustration},  for $ p $-adaptive flux reconstruction methods, we nest polynomial degrees on $ p $-sublevels towards the lower ones for a given order gap between different sublevels. For example, for a $ p $-adaptive $ p^4 $ FR method, when the order gap 2 is used for 3 $ p $-sublevels, the polynomial hierarchies for elements on the finest level which has $ \mathbb{P}^4 $, $ \mathbb{P}^3 $, $ \mathbb{P}^2 $, and $ \mathbb{P}^1 $ will be $ p\{4\text{-}2\text{-}0\} $, $ p\{3\text{-}1\text{-}0\} $, $ p\{2\text{-}1\text{-}0\} $, and $ p\{1\text{-}0\text{-}0\} $, respectively; when the order gap is 1,  the polynomial hierarchies for elements on the finest level which has $ \mathbb{P}^4 $, $ \mathbb{P}^3 $, $ \mathbb{P}^2 $, and $ \mathbb{P}^1 $ will be $ p\{4\text{-}3\text{-}2\} $, $ p\{3\text{-}2\text{-}1\} $, $ p\{2\text{-}1\text{-}0\} $, and $ p\{1\text{-}0\text{-}0\} $, respectively. 
    With this being said, when it comes to the polynomial hierarchy of $ p $-adaptive FR methods, we will only provide that of the elements in the flow field which has polynomial degree $ p_\mathrm{max} $.
\section{Numerical experiments}
    In this section, we present the results from numerical investigations of $p$MG methods, first when applied as a nonlinear solver in steady problems, then when applied as a precondition for unsteady problems. The $ p $-adaptation methods of \cref{coupling_adaptation} were coupled to $p$MG for both steady and unsteady problems. Furthermore, the smoothing used on the finest $p$-sublevel was exclusively the element-Jacobi method; however, the smoother at the coarsest $p$-sublevel was subject to investigation.
    
    Throughout these investigations, the flow field was initialized with the impulsive condition of $ (\rho, u, v, \text{Ma})^\top = (1,1,0,\text{Ma}_\infty)^\top $. All numerical experiments were conducted system with an Intel Core i7-10750H CPU @ \SI{2.6}{\giga\hertz}, and $ 2\times8 $GB DDR4 RAM @ \SI{3200}{\mega\hertz}. All simulations were run using four processes.
    
\subsection{$p$MG for steady problems}\label{sec:num_steady}

\subsubsection{Inviscid flow over a sphere}
\begin{figure}[tbhp]
    	\centering
    	\subfloat[Normalised pressure.]{\label{fig:sphere_contour_p}	\adjustbox{width=0.48\linewidth, valign=b}{\includegraphics[width=1\linewidth]{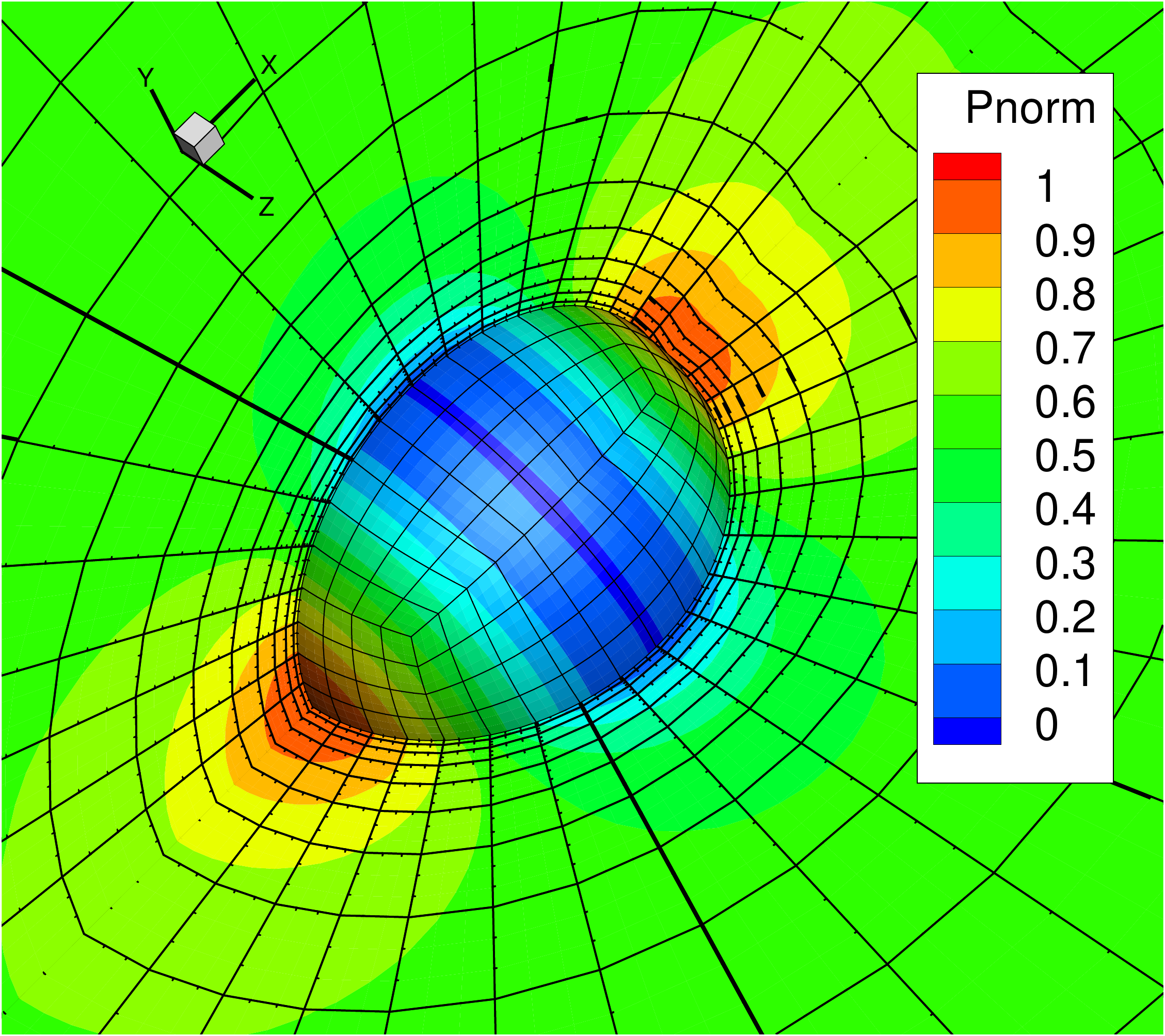}}}
    	~
    	\subfloat[Mach number.]{\label{fig:sphere_contour_ma}	\adjustbox{width=0.48\linewidth, valign=b}{\includegraphics[width=1\linewidth]{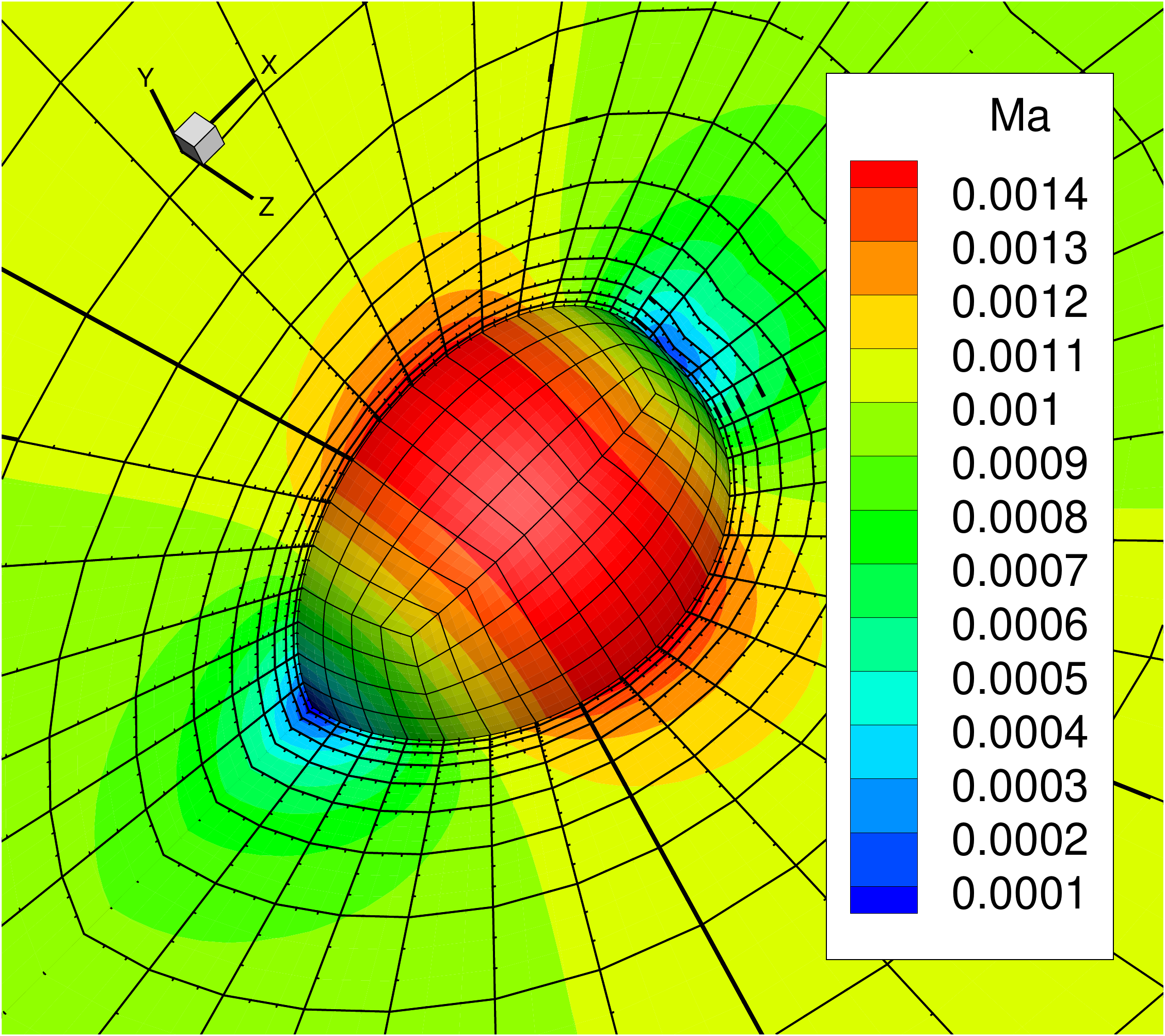}}}
    	\caption{\label{fig:sphere_contour}Contoured quantities for the inviscid flow past a sphere at $\mathrm{Ma}_\infty=0.001$.}
    \end{figure}

    We first investigated the performance of $p$MG methods for steady problems by simulating the inviscid flow over a sphere at low Mach number, where $\mathrm{Ma}_\infty=0.001$. The diameter of the sphere was set as $d=1$, with only a quarter of the sphere being meshed with hexahedral elements. Symmetric boundary conditions were applied to the symmetric planes and far field boundary conditions were applied to the outer boundaries. Inviscid/slip wall boundary conditions were imposed on the wall. A mesh with a total of \num{2464} quadratic hexahedral elements  was used, and the near-wall mesh used for this simulation is shown in \cref{fig:sphere_contour} as well as contours of the normalized pressure and $ \mathrm{Ma} $. To preserve accuracy and accelerate convergence, the local preconditioning technique originally developed by Weiss and Smith~\cite{weiss1995preconditioning} was applied~\cite{wang2019implicit}. At the time of writing, the authors have not observed any acceleration by applying this local preconditioning to JFNK methods. Therefore, local preconditioning is restricted to only steady problems.
    
    In tackling this problem, $p$-adaptive and $p$-uniform methodologies were applied to $\mathbb{P}^3$ and $\mathbb{P}^4$ FR to investigate how the smoothing method at the coarsest $p$-sublevel affects the convergence of $p$MG. A 3-level $p$MG nonlinear solver was used with an order gap of two between adjacent $p$-sublevels. Specifically, the polynomial hierarchy on $\mathbb{P}^3 $ and $\mathbb{P}^4$ elements was $ p\{3\text{-}1\text{-}0\} $ and $ p\{4\text{-}2\text{-}0\} $, respectively. The number of smoothing iterations at each level was $n\{2\text{-}4\text{-}5\}$. Note that the predicted drag coefficient for flows past a sphere is defined as
    \begin{equation}\label{eq:drag_coeff}
     C_d = \frac{2F_d}{\rho_\infty U_\infty^2A},\quad  \text{and}\quad A = \pi d^2/4.
    \end{equation} 
    
\begin{figure}[thbp]
		\centering
		\subfloat[]{\label{fig:sphere_degree_p3}\adjustbox{width=0.48\linewidth,valign=b}{\includegraphics[width=1\linewidth]{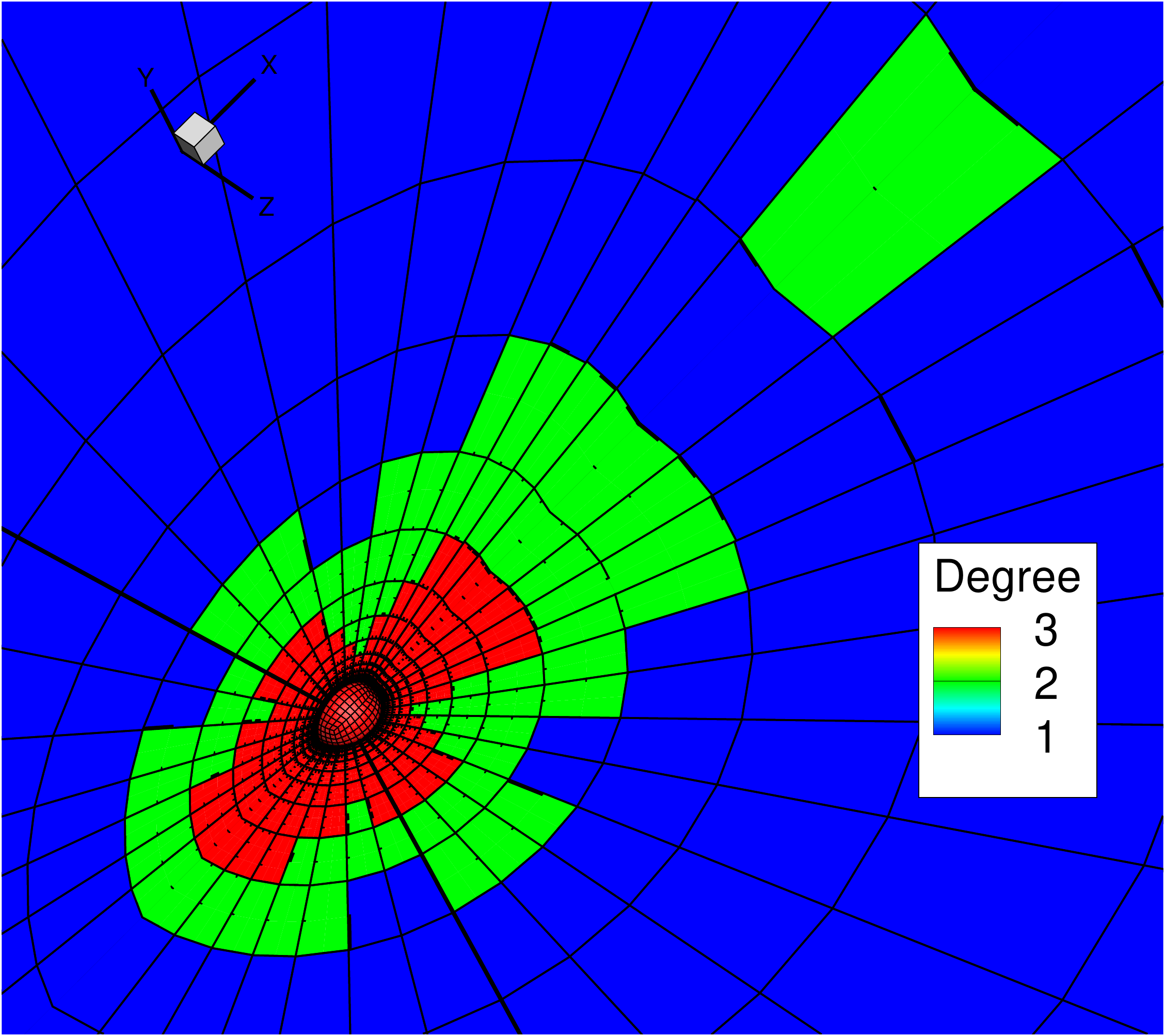}}}
		\quad
		\subfloat[]{\label{fig:sphere_degree_p4}\adjustbox{width=0.48\linewidth,valign=b}{\includegraphics[width=1\linewidth]{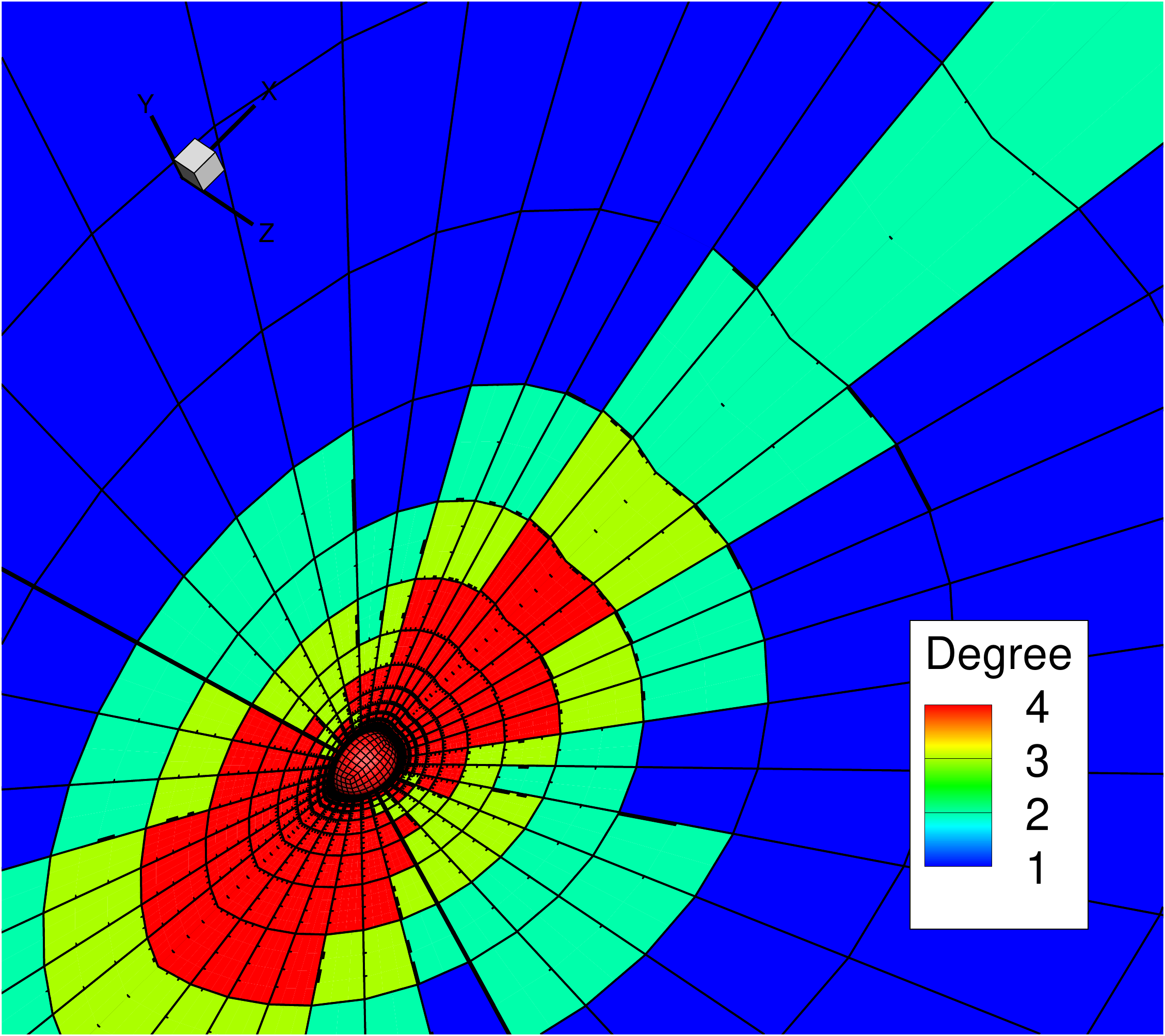}}}
		\caption{\label{fig:p-distribution-sphere}Distribution of polynomial degrees of $ p $-adaptive (a) $ \mathbb{P}^3 $ and (b) $ \mathbb{P}^4 $ FR  solving inviscid flow over sphere at $ \mathrm{Ma}_\infty=0.001 $.}
	\end{figure}

    The results presented in \cref{fig:p-distribution-sphere} show no pressure oscillations near the front stagnation point, which demonstrates that the modified Roe solver in the local preconditioning technique~\cite{weiss1995preconditioning} is able to balance the dissipation on different characteristics of the convection terms. Furthermore, the polynomial distributions presented in \cref{fig:p-distribution-sphere} show that, upon convergence, the  $p$-adaptive field is symmetric for both $\mathbb{P}^3$ and $\mathbb{P}^4$ FR. This is consistent with the symmetric solution. The predicted drag coefficients obtained via $p$-adaptive $\mathbb{P}^3$ and $\mathbb{P}^4$ FR  are \num{1.25e-5} and \num{5.01e-6}, respectively. Those obtained from $ p $-uniform $\mathbb{P}^3 $ and $\mathbb{P}^4 $ FR are \num{4.42e-7} and \num{1e-7}. 
    
    It has previously been shown that feature-based adaptation methodologies are inferior to output-based ones~\cite{li2011continuous, fidkowski2011review} for steady problems. This investigation is only intending to study the behaviour of $p$MG when coupled to $p$-adaptation, and furthermore, the simple feature-based $p$-adaptation approach is advantageous for unsteady problems in terms of implementation. Further discussion of $p$-adaptation methodologies is beyond the scope of this study.
    
\begin{figure}[tbhp]
	    \centering
	    \subfloat[Absolute residual vs. number of V-cycles.]{\label{fig:sphere_degree_p3_res_vcycle}\adjustbox{width=0.48\linewidth,valign=b}{\includegraphics{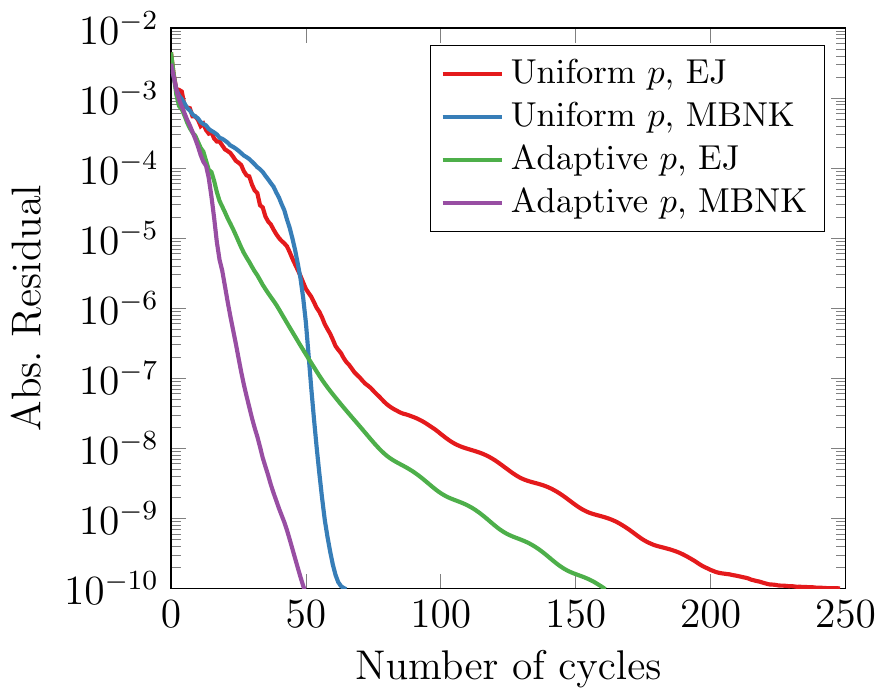}}}
	    \quad
		\subfloat[Absolute residual vs. wall time.]{\label{fig:sphere_degree_p3_res_cpu}\adjustbox{width=0.48\linewidth,valign=b}{\includegraphics{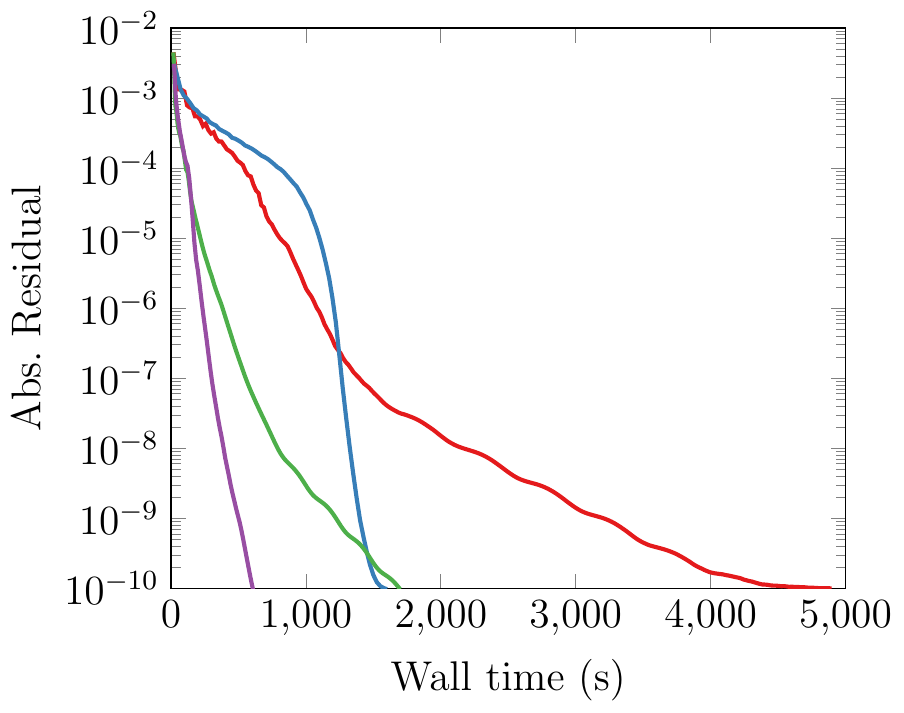}}}
        \caption{\label{fig:sphere_p3_residual}Solver statistics for inviscid flow over sphere at $ \mathrm{Ma}_\infty=0.001$ using $\mathbb{P}^3$ FR in several configurations.}
	\end{figure}

\begin{figure}[tbhp]
	    \centering
	    \subfloat[Absolute residual vs. number of V-cycles.]{\label{fig:sphere_degree_p4_res_vcycle}\adjustbox{width=0.48\linewidth,valign=b}{\includegraphics{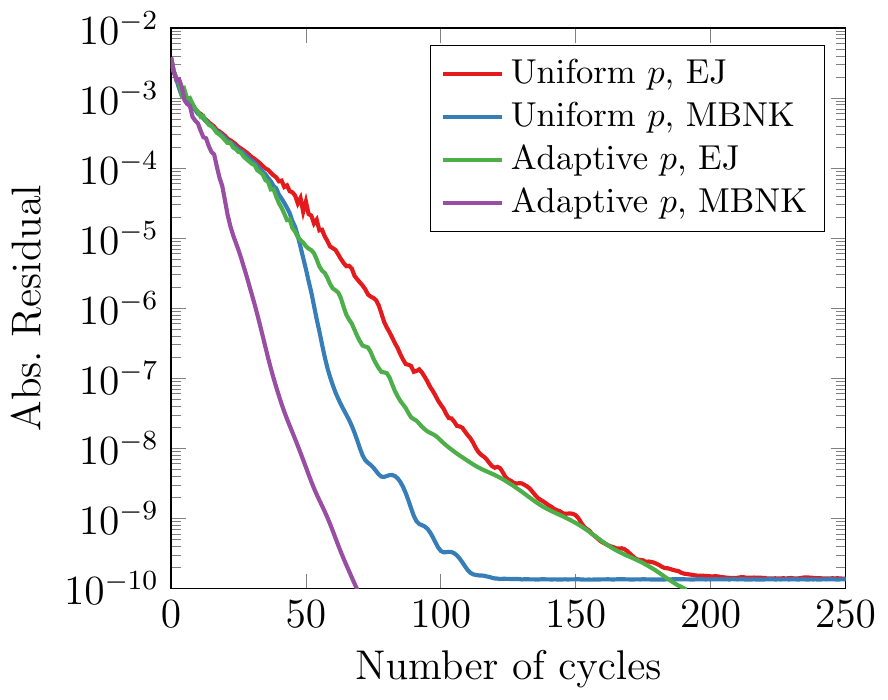}}}
	    \quad
		\subfloat[Absolute residual vs. wall time.]{\label{fig:sphere_degree_p4_res_cpu}\adjustbox{width=0.48\linewidth,valign=b}{\includegraphics{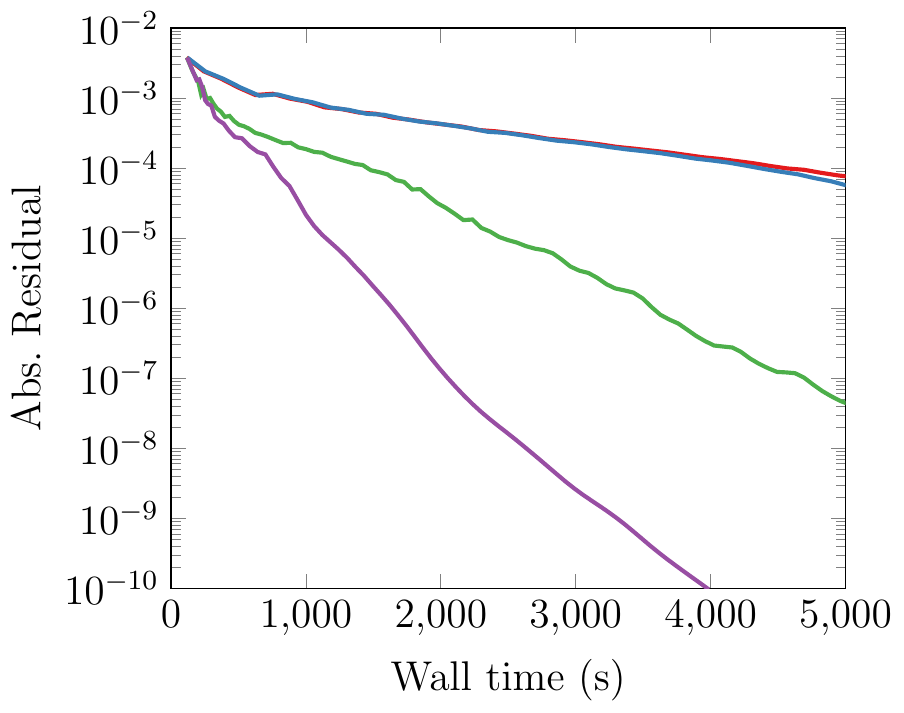}}}
        \caption{\label{fig:sphere_p4_residual}Solver statistics for inviscid flow over sphere at $ \mathrm{Ma}_\infty=0.001$ using $\mathbb{P}^4$ FR in several configurations.}
	\end{figure}

    Statistics of $p$-uniform and $p$-adaptive FR methods with EJ and MBNK smoothing at the coarsest level are documented in \cref{tab:inv_sphere_stats}. From the residual data presented in \cref{fig:sphere_p3_residual,fig:sphere_p4_residual}, it is clear that the residual in some calculations stalled at ${\sim}\num{1e-10}$. For this reason the wall clock time data present in \cref{tab:inv_sphere_stats} is the time taken for the absolute residual to drop below \num{5e-10}. For both $p$-adaptive and $p$-uniform $\mathbb{P}^3 $ FR, a speedup of ${\sim}2.5$ was observed when using the MBNK smoother rather than the EJ smoother on the coarsest $p$-sublevel in $p$MG. Similarly for $\mathbb{P}^4$ FR, MBNK achieved a speedup of  $1.59$ and  $2.79$ for $p$-uniform and $p$-adaptation, respectively. Overall, the significance of employing a matrix-based smoother on the coarsest $p$-sublevel is consistent for $p$-uniform and $p$-adaptive FR methods.
    
    \begin{table}[tbhp]
        \centering
        \caption{Statistics of $p$-adaptive and $p$-uniform FR solving invscid flow over the sphere at $\text{Ma}_\infty=0.001$ for different smoothers at the coarsest level. Speed-up is shown for MBNK relative to EJ.}
        \label{tab:inv_sphere_stats}
        \begin{tabular}{c c r r c}
    	    \toprule
            $p$ & Adaptation  & \multicolumn{2}{c}{Runtime (s)} & Speedup\\
            & & EJ & MBNK & \\
            \midrule
            $3$ &No  & 3462 & 1451  & 2.49 \\
            $3 $ & Yes & 1364 & 549 & 2.48\\
            $4 $ & No & 20053 & 12622 & 1.59\\
            $4 $ & Yes & 9919 & 3549 & 2.79\\
            \bottomrule
        \end{tabular}
    \end{table}
    
	In the work of Shahbazi et al.~\cite{shahbazi2009multigrid}, a more efficient smoother at the coarsest $p$-sublevel (either $\mathbb{P}^0 $ or $\mathbb{P}^1 $) was produced by building up deeper $hp$MG hierarchy and employing saw-tooth/W-cycle. They found sufficient smoothing on $\mathbb{P}^0$ is critically important when solving inviscid problems while sufficient smoothing on $\mathbb{P}^1$ is more important for solving viscous problems, and the results presented in \cref{fig:sphere_p3_residual,fig:sphere_p4_residual} are consistent with this conclusion. Naturally, this benefit motivated us to further improve the smoother quality at intermediate $p$-sublevel when the convergence rate of PTC stalls for viscous problems which was observed in~\cite{wang2019p}. 
	
\subsubsection{Viscous flow over a sphere}	
	Using the same mesh as shown in \cref{fig:sphere_contour}, we considered the case of $\mathrm{Re}_d=118 $ and $\text{Ma}_\infty=0.001$. Three different strategies were used for $\mathbb{P}^3$ FR with $p$-adaptation:
	\begin{description}
	    \item[Strategy 1.] The first strategy is 3-level $p$MG cycle with an order gap 1, where $\mathbb{P}^3$ elements have a $p \{3\text{-}2\text{-}1\}$ polynomial hierarchy. EJ is used on all levels except the coarsest level where MBNK is used;
	    \item[Strategy 2.]  The second one uses the same polynomial hierarchy and smoother; however, after 15 pseudo-time iterations the smoother on the intermediate $p$-sublevel will be switched from EJ to MBNK; 
	    \item[Strategy 3.] The third strategy employs a 2-level $p$MG cycle only, which for $\mathbb{P}^3$ elements is a $p\{3\text{-}2\}$ polynomial hierarchy. Similarly, EJ smoothing is used on the finest level, and MBNK on coarsest level.
	\end{description}
    Regarding the number of iterations for smoothing,  $ n\{10\text{-}10\} $ and $ n\{10\text{-}10\text{-}10\} $ are used for 2-level and 3-level $p$MG, respectively.
    
\begin{figure}[tbhp]
	    \centering
	    \subfloat[Absolute residual vs. number of V-cycles.]{\label{fig:sph_p3_vcycle}\adjustbox{width=0.48\linewidth,valign=b}{\includegraphics{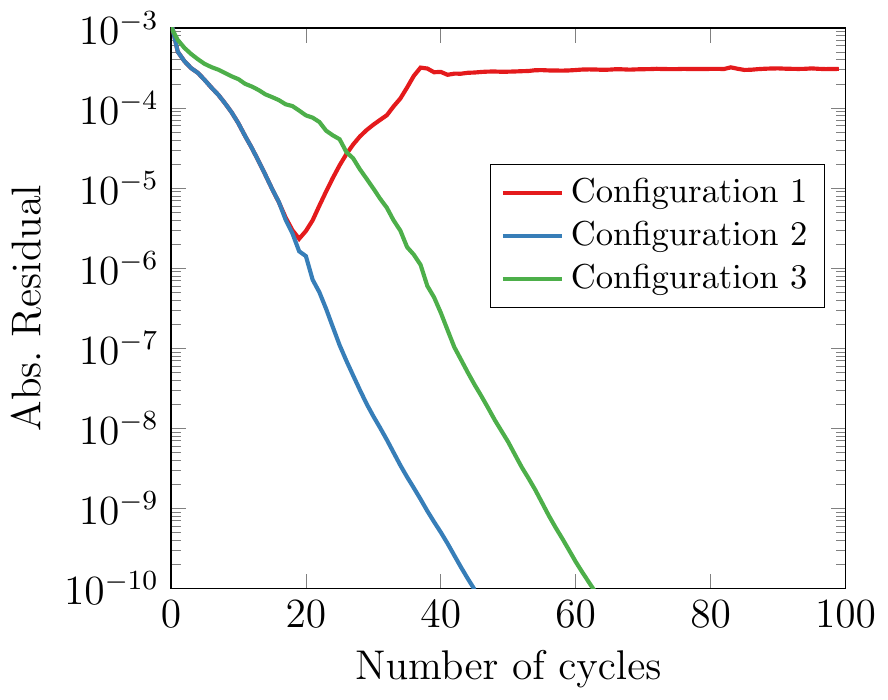}}}
	    ~
	    \subfloat[Absolute residual vs. wall time.]{\label{fig:sph_p3_time}\adjustbox{width=0.48\linewidth,valign=b}{\includegraphics{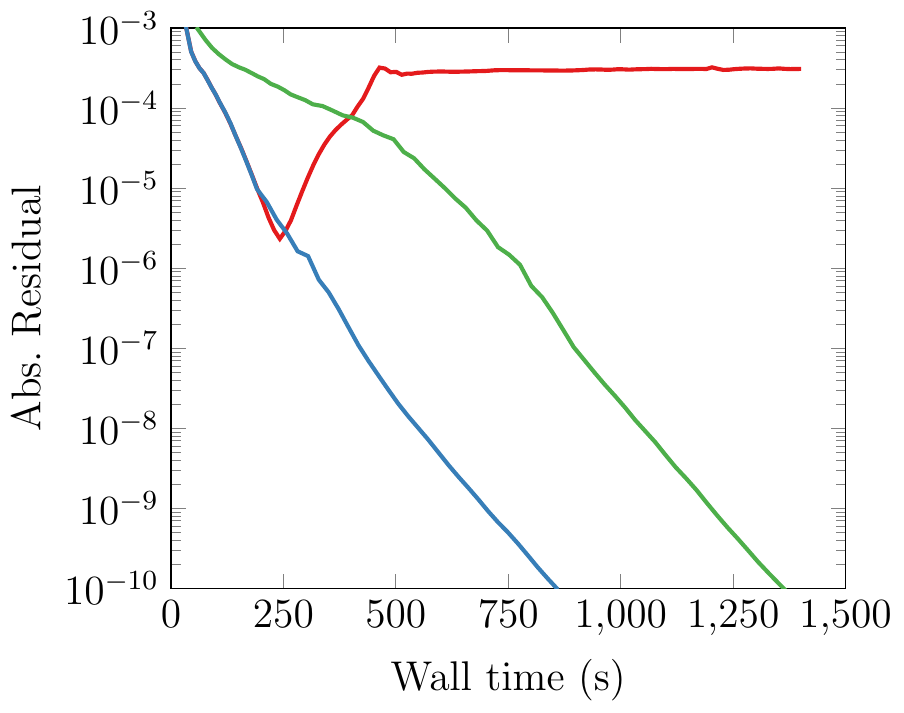}}}
	    \\
	    \subfloat[$\Delta\tau$ vs. number of V-cycles.]{\label{fig:sph_p3_dtau}\adjustbox{width=0.48\linewidth,valign=b}{\includegraphics{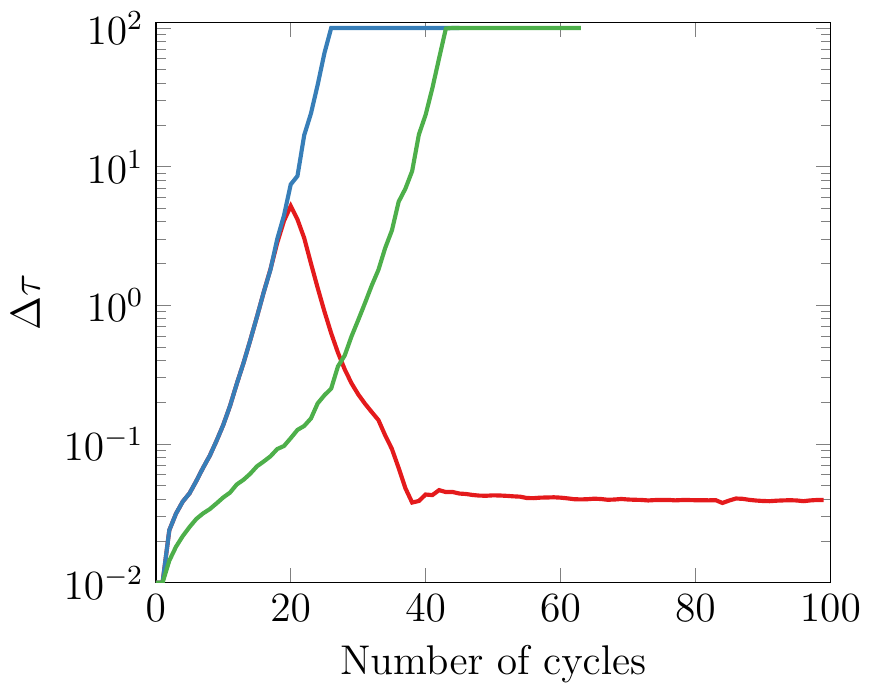}}}
	    ~
	    \caption{Effect of $p$MG smoother configuration for viscous flow over sphere at $ \text{Ma}_\infty=0.001$ and $\mathrm{Re}_d=118$ using $p$-adaptive $\mathbb{P}^3$ FR.}
    	\label{fig:sphere_viscous}
	\end{figure}

    The predicted drag coefficient is $0.9700$ which is slightly smaller than 1.002 obtained in~\cite{ji2021p} at $ \text{Ma}=0.2535 $. The residual histories of these three solution strategies for $ p $-adaptive $ p^3 $ FR are depicted in Figure~\ref{fig:sphere_viscous}. Ideally, as residual drops, pseudo time step $ \Delta \tau $ will keep increasing. However, for the first strategy, when $ \Delta \tau $ got close to $ 10 $, some instability was triggered and the  residual hence increased to a large value ($ \sim 10^{-3} $) and stalled at such level. For the second strategy, we switched the smoother on $\mathbb{P}^2$ sublevel from EJ to MBNK after 15 V-cycles or pseudo time iterations.  With better smoothing on intermediate $p$-sublevel, PTC was stabilized. It is  unclear how such switch would  quantitatively affect the CFL condition. When comparing  {Strategy 3} with  {Strategy 2}, surprisingly the latter was faster in terms of both convergence rate and convergence speed. We hypothesise that in PTC, at the initial stage the errors that dictates the convergence rate are those low-frequency modes. As the residual drops, the frequency of dominating ones gradually gets higher. Consistent with our assumption,  {Strategy 3} with MBNK on $ \mathbb{P}^2 $ actually had smaller convergence rate initially. This simple numerical experiments raises up two questions to be further explored, (a) how smoothing on sublevels would affect the stable CFL of $p$MG and (b) how to design the smoothing strategy such as when to switch the smoother on intermediate sublevels.
    
\clearpage
    
\subsection{$p$MG preconditioner for unsteady problems}\label{sec:num_unsteady} 
    In this section we focus on stiff systems arising from low-Mach-number flows and high-aspect-ratio elements in near wall regions. We present results from investigations which aimed to demonstrate that JFNK can be greatly accelerated for these problem by using preconditioners other than EJ, in particular an advanced $ p $MG preconditioner. In this subsection, we used MNBK as the coarsest sublevel smoother for all $p$MG solvers and preconditioners.
    
\subsubsection{Unsteady flow over a cylinder}
    The first case was 2D unsteady flow over a circular cylinder at $ \mathrm{Re}_d=\num{1200}$ when $\mathrm{Ma}_\infty=0.01 $, for diameter $d = 1$. A similar problem was studied in~\cite{bijl2002implicit, wang2020comparison} to evaluate the accuracy and efficiency of implicit time integration methods at higher $\mathrm{Ma}$ numbers. 
    
    The mesh used in this investigation had \num{5316} quadratic quadrilateral elements, with the first layer of near-wall elements having a height of 0.01. This gave an \textit{a prior} estimation of $y^+\approx0.8$. For high-order numerical methods, a different opinion in terms of the \textit{a prior} estimation of the element size is that one could use the distance of the first layer solution points such that the element size can be greatly increased to give a posterior estimation of $ y^+ $ close to 1 at these points. However, throughout this work, we followed the convention for low-order methods.
    We aimed to investigate the efficacy of $p$MG methods for stiff systems in the low Mach number regime. We used ESDIRK2 with a $p$MG preconditioner to run the simulation for 150 convective time units and employed the solutions as initial values for our numerical study. Note that no $p$-adaptation was used for simulations in this subsection. 
    
\begin{figure}[tbhp]
    	\centering
    	\adjustbox{width=0.5\linewidth,valign=b}{\includegraphics{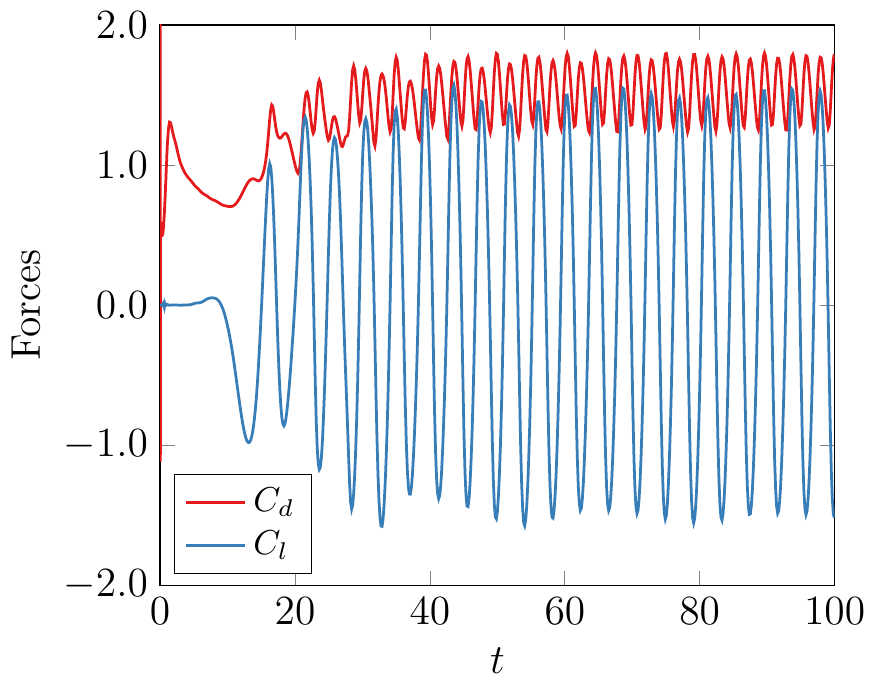}}
    	\caption{\label{fig:naca_forces}Lift and drag coefficient for the viscous flow over 2D circular cylinder at $\mathrm{Re}_d=\num{1200}$ and $ \text{Ma}_\infty=0.01 $, using $p$-uniform $\mathbb{P}^3$ FR.}
    \end{figure}

\begin{figure}[tbhp]
    	\centering
    	\begin{subfigure}[]{1\linewidth}
    		\centering
    		\includegraphics[width=0.8\linewidth]{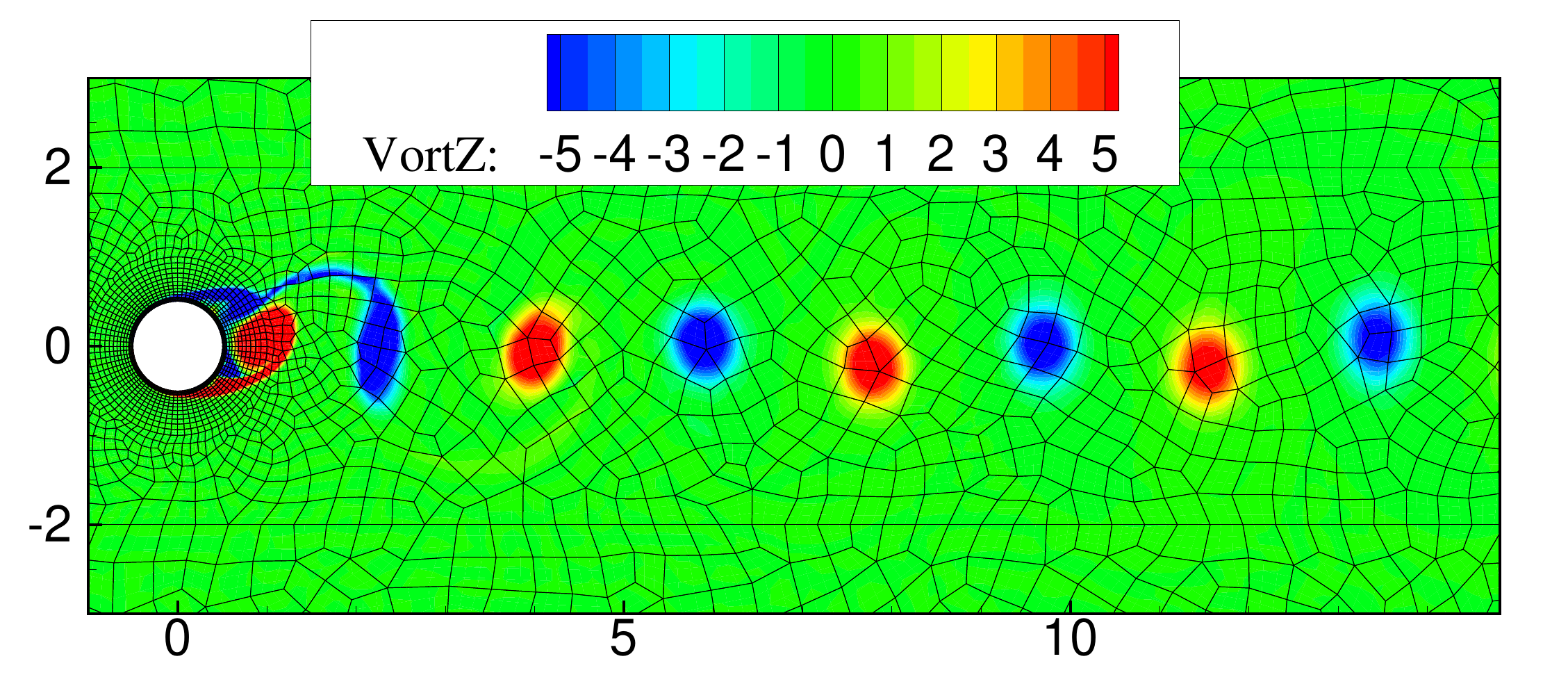}
    		\caption{\label{fig:naca_a}z-Vorticity.}
    	\end{subfigure}
    	
    	\begin{subfigure}[]{1\linewidth}
    		\centering
    		\includegraphics[width=0.8\linewidth]{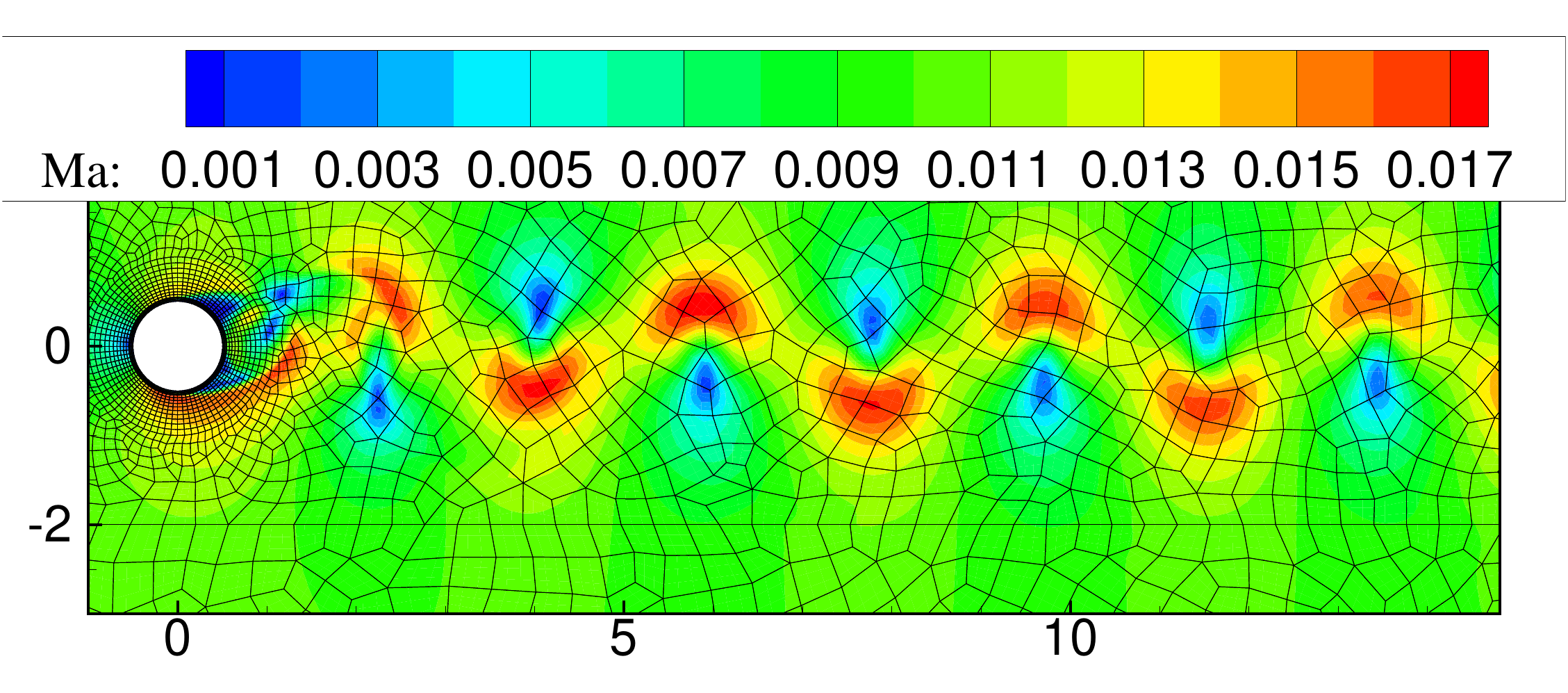}
    		\caption{\label{fig:naca_b}$\mathrm{Ma}$}
    	\end{subfigure}
    	\caption{Contours of instantaneous quantities for the viscous flow over a 2D circular cylinder at $\mathrm{Re}_d=\num{1200}$ and $\mathrm{Ma}_\infty=0.01$.}
    	\label{fig:naca}
    \end{figure}

\begin{figure}[tbhp]
    	\centering
    	\includegraphics[width=0.7\linewidth]{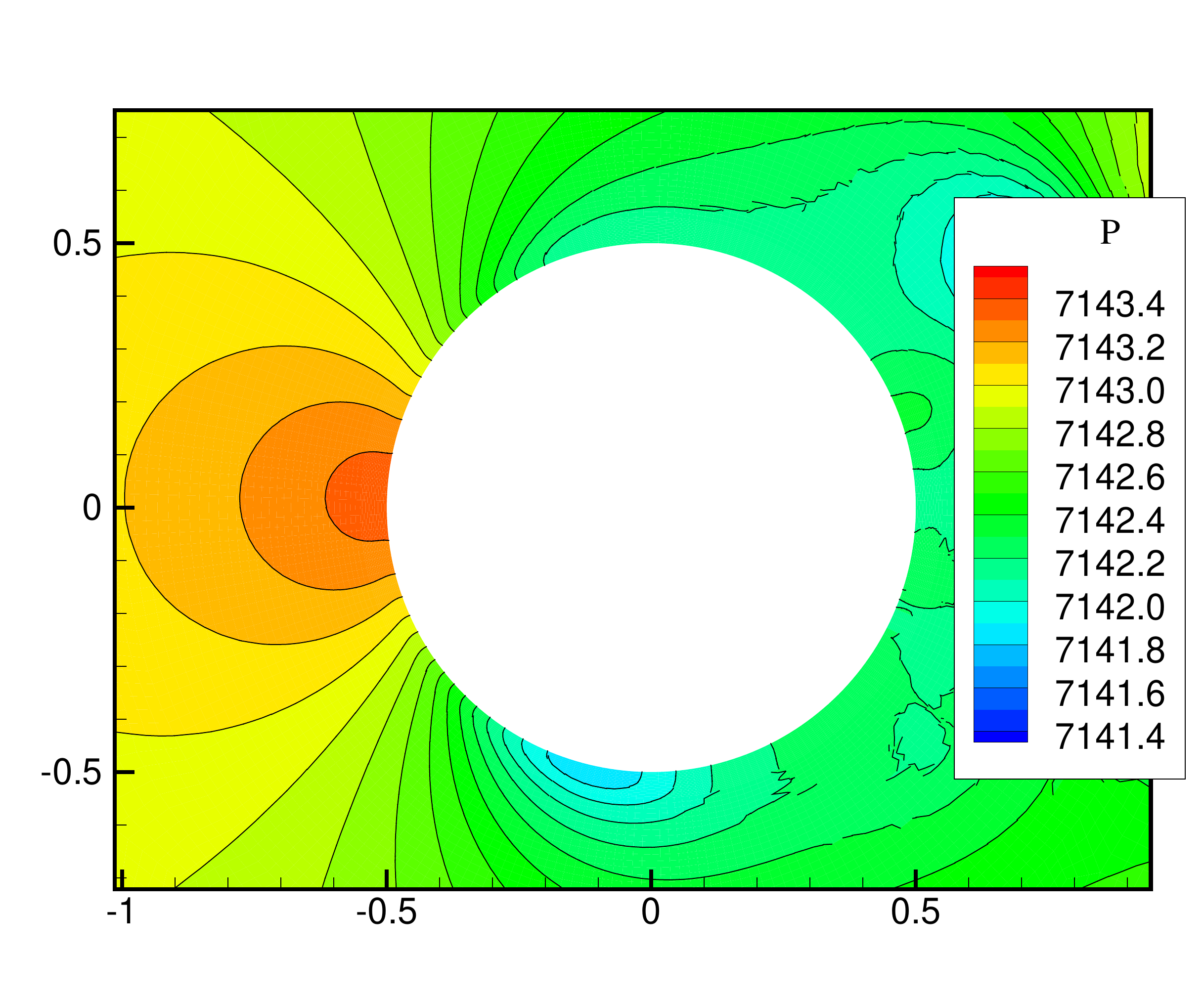}
    	\caption{Pressure field near the stagnation point of  viscous flow over 2D circular cylinder at $\mathrm{Re}_d=\num{1200}$ and $ \mathrm{Ma}_\infty=0.01$.}
    	\label{fig:naca_pressure}
	\end{figure}

    The predicted lift and drag coefficients are illustrated in~\cref{fig:naca_forces}. The predicted Strouhal  number, defined as ${fd}/{U_\infty}$, is 0.2386, which is slightly smaller than the predicted value 0.2467 in~\cite{bijl2002implicit} at $\mathrm{Ma} = 0.3$.  Instantaneous contours of vorticity and $\mathrm{Ma}$ are shown in \cref{fig:naca} as well as the local mesh. In \cref{fig:naca_pressure}, near the front stagnation point, we have not observed any pressure oscillations, even though only the original Roe Riemann solver~\cite{RoeSolver} was employed. For all following comparison in this subsection, we  documented  the statistics for 5 convective time units.
    
\clearpage
    Three solution strategies were explored first: $p$MG as a nonlinear solver with the cycles $ p\{3\text{-}1\text{-}0\} $ and $ p\{3\text{-}2\text{-}1\text{-}0\} $ for the ESDIRK2 scheme, as well as  JFNK with EJ preconditioning as a comparison. The $p$MG nonlinear solver used smoothing steps $ n\{2\text{-}2\text{-}5\}$ and $n\{2\text{-}2\text{-}2\text{-}5\}$ for 3-level and 4-level cycles respectively. A limit on pseudo-time iterations was set as $N_\mathrm{ptc}^\mathrm{max}= 100$ and relative tolerance of PTC,  $\mathrm{tol}_\mathrm{ptc}^{r}=10^{-4}$, was used. When using $p$MG as a nonlinear solver in this case, the smoothers were updated every 20 pseudo time iterations due to the rapid decrease in the convergence rate observed in this case.
    
    For JFNK-EJ, we set the relative tolerance for GMRES as $ \mathrm{tol}_\mathrm{gmres}^\mathrm{r}=10^{-2} $. A larger one, such as $ \mathrm{tol}_\mathrm{gmres}^\mathrm{r}=10^{-1} $, would lead to instability of pseudo transient continuation for this simple EJ preconditioner. We only updated the smoothers each physical time step and the convergence of PTC was remarkable.  $ K_\mathrm{dim} = 30 $ was set as the maximum dimension of the Krylov subspace to make sure that GMRES would not diverge. 
    
\begin{figure}[tbhp]
    	\centering
        \subfloat[Relative Residual]{\label{fig:stage2_naca_a}\adjustbox{width=0.48\linewidth,valign=b}{\includegraphics{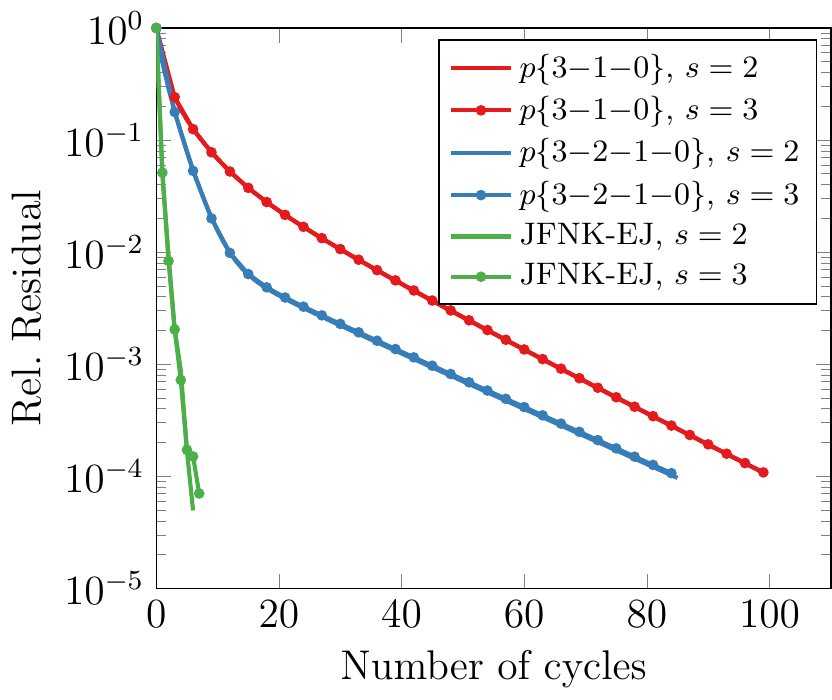}}}
        \quad
        \subfloat[Absolute Residual]{\label{fig:stage2_naca_b}\adjustbox{width=0.48\linewidth,valign=b}{\includegraphics{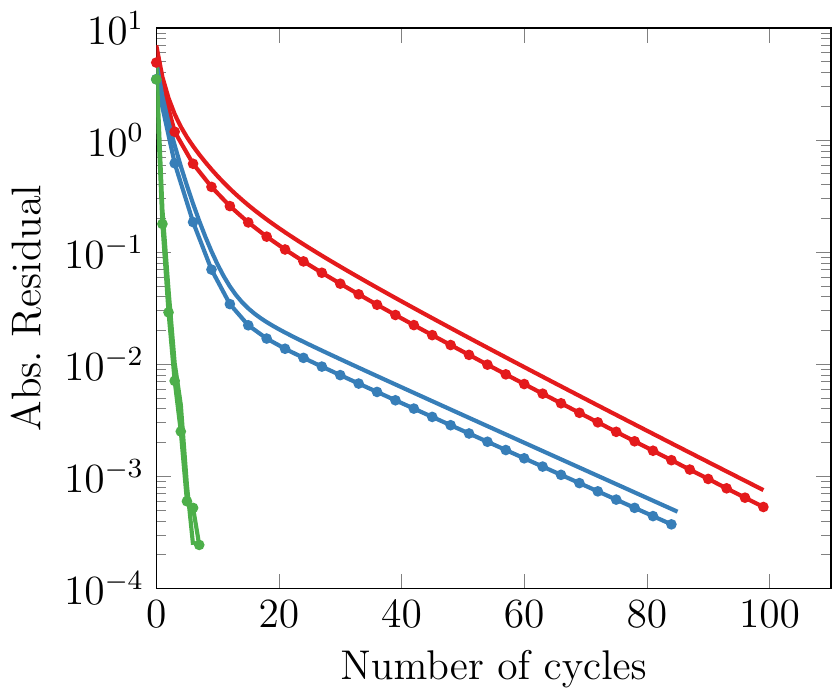}}}
    	\caption{Residual evolution for the 2D viscous flow over a 
    	cylinder at $\mathrm{Re}_d=\num{1200}$ and $\mathrm{Ma}_\infty=0.01$, using $p$MG and JFNK nonlinear solver, with residuals displayed for different ESDIRK2 stages.}
    	\label{fig:stage2_naca}
    \end{figure}

    Statistics of numerical simulations using the $ p $MG nonlinear solver and JFNK with the EJ preconditioner are documented in \cref{tab:pmg_vs_pmg_precon}.
    This shows a small runtime benefit to the $p$MG nonlinear solver with $p\{3\text{-}1\text{-}0\}$ over that with $p\{3\text{-}2\text{-}1\text{-}0\}$. Presented in \cref{fig:stage2_naca} are typical residual histories for all three methods. This shows that even though $p\{3\text{-}2\text{-}1\text{-}0\}$ led to a larger initial convergence rate, the overhead of the extra $ \mathbb{P}^2 $ sublevel indeed made the $p$MG nonlinear solver slightly slower. This is consistent with observations from the previous steady problem where Strategy 3 was slower than Strategy 2. Moreover, JFNK-EJ was 2.26 times faster than the baseline case $ p $MG with $ p\{3\text{-}1\text{-}0\} $. A further observation was that $ p $MG nonlinear solvers needed an order of magnitude more pseudo time iterations to achieve the same residual drop of 4 orders of magnitude. On the contrary, with JFNK-EJ the convergence rate of PTC was reminiscent of Newton-like approaches and pseudo transient continuation converged on average within 8 pseudo time iterations for each stage of ESDIRK2 ($ N_\mathrm{ptc}^\mathrm{avg} = 7.3 $); meanwhile the average number of GMRES iteration for each stage, $N_\mathrm{gmres}^\mathrm{avg}$, was significant. For both polynomial hierarchies investigated, the convergence rate for the $ p $MG nonlinear solver got smaller after the relative residual dropped lower than $10^{-2}$. For unsteady flow simulations, the initial guess of the nonlinear system --- which is actually the solution at time step $t^n$ --- is close to the solution at $t^{n+1}$. That is to say, the nonlinear solver started at the stage of PTC when errors on intermediate $ p $-sublevel dictated convergence. Therefore, extra smoothing on the $ \mathbb{P}^2 $ sublevel gave an improved convergence rate initially. However, as the residual dropped, errors on the $ \mathbb{P}^3 $ sublevel began to dominate and the overhead of the $\mathbb{P}^2$-sublevel was largely wasted. These results shows that tuning the configuration of the $ p $MG nonlinear solver is unlikely to make it outperform JFNK-EJ.

    We further tested how a $ p $MG preconditioner can improve the performance of JFNK used in ESDIRK. The time step was increased from $ \Delta t=0.025 $ to $ \Delta t=0.05 $. With $ \Delta t $ being doubled, the stiffness of the nonlinear system at each stage of ESDIRK was increased due to larger source terms. For this larger times steps the lower tolerance of $\mathrm{tol}_\mathrm{gmres}^\mathrm{r}=10^{-3}$ was used to stabilize JFNK-EJ with a significantly larger dimension of $ K_\mathrm{dim}=100 $. The $ p $MG preconditioner coupled to JFNK was $ p\{3\text{-}1\} $ with $ n\{2\text{-}2\}$. The tolerance and dimension of GMRES for JFNK with a smaller time step were set as $ \mathrm{tol}_\mathrm{gmres}^\mathrm{r}=10^{-1} $ and $K_\mathrm{dim}=30$ with no instability being observed in PTC. Statistics of this test are documented in \cref{tab:ESDIRK-JFNK-EJ-PMG}. For JFNK-EJ, increasing the time step decreased its efficiency. It is observed from \cref{tab:ESDIRK-JFNK-EJ-PMG} that, even though the averaged number of pseudo time iterations of PTC per RK stage $N_\mathrm{ptc}^\mathrm{avg}$ was only slightly increased, the averaged number of GMRES iterations per RK stage $ N_\mathrm{gmres}^\mathrm{avg} $ was increased by a factor greater than 2. On the contrary, for the larger time step $ \Delta t =0.05$ JFNK-$p$MG was almost 2 times faster than JFNK-EJ and more than 4.5 times faster than the baseline case; only 44 GMRES iterations were needed per RK stage. Moreover, the dimension of the Krylov subspace was maintained as that of JFNK-EJ when $ \Delta t=0.025 $.    
    
    Then moving on to consider the linearly implicit ROW scheme, as suggested in the study~\cite{wang2020comparison}, GMRES has to converge the residual to a relatively smaller value. With this in mind, the GMRES tolerence of $ \mathrm{tol}_\mathrm{gmres}^\mathrm{r}=10^{-6} $ was used when investigating the benefit of $ p $MG preconditoners. The simulation was run at $ \Delta t = 0.05 $ when using ROW2  with both the EJ and $p$MG preconditioner, with different configurations of polynomial hierarchy was also investigated.  Statistics of these tests are documented in \cref{tab:ROW-JFNK-EJ-PMG}. 
    A significant result is that the advanced $ p $MG preconditoner with $p$-hierarchy $p\{3\text{-}1\}$ can achieve a speedup of  5.32 compared to the baseline and 1.97 when compared to the EJ preconditioner. 
    $ N_\mathrm{gmres}^\mathrm{avg} $ per RK stage was below 31, one order of magnitude smaller than that of GMRES-EJ. This data also shows how a bad choice of $ p $-hierarchy can deteriorate the performance of $ p $MG preconditioner. For example, from \cref{tab:ROW-JFNK-EJ-PMG} GMRES-$ p $MG with $ p\{3\text{-}0\} $ not only required a large $ K_\mathrm{dim} =200 $ for stability, but also required greater than three times the number of GMRES iterations per RK stage to converge. Furthermore, the performance of this strategy was worse than that of the EJ preconditioner. This observation is similar to that of \cref{sec:num_steady} where insufficient smoothing of the dominant error modes would adversely effect the convergence. A final observation was that $N_\mathrm{gmres}^\mathrm{avg}$ of GMRES-$ p $MG with $ p\{3\text{-}2\text{-}1\} $ was notably smaller than that of GMRES-$ p $MG with $p \{3\text{-}1\} $. This could be expected since the additional smoothing on $\mathbb{P}^2$ better addresses the error modes on this intermediate sublevel. However, the overhead on this sublevel slightly decreased its efficiency.

     To summarize, for a stiff system in the low-Mach-number regime, JFNK is more preferable than $ p $MG to solve the nonlinear system arising from implicit time stepping methods. ESDIRK2 and ROW2 can be 4.56 and 5.32 times faster, respectively, when GMRES-$ p $MG is used instead of $ p $MG nonlinear solver. The developed $ p $MG preconditioner can significantly speed up the simulation by factor of 2.01 and 1.97 for ESDIRK2 and ROW2, respectively, when compared to EJ preconditioner with proper configuration of the polynomial hierarchy. Lower  RAM usage of $ p $MG preconditioner is realized via using one order of magnitude smaller dimension number for the Krylov subspace. Similar to the steady problems considered in this work, insufficient smoothing on intermediate $ p $-sublevels will actually worsen the performance of $ p $MG preconditioners greatly. Even though a dense $ p $MG can mitigate the problem, the overhead on relatively high $ p $-sublevel would comprise the efficiency. We genuinely would recommend a $ p\{p_0\text{-}(p_0/2)\} $ polynomial hierarchy for this type of problem. With a better preconditioner, $ p $MG in our case, ROW2 was about 1.16 times faster than ESDIRK2.
    
    \begin{table}[tbhp]
    	\centering
    	\caption{$ p $MG nonlinear solvers vs. JFNK-EJ for ESDIRK2 for viscous flow over a cylinder at $ \text{Ma}_{\infty} =0.01 $. $ N_\mathrm{ptc}^\mathrm{avg} $ is the averaged number of pseudo time iterations per RK stage and $ N_\mathrm{gmes}^\mathrm{avg} $ is the averaged number of GMRES iterations per RK stage.}
    	\label{tab:pmg_vs_pmg_precon}
    	\begin{tabular}{c c c c c c c}
    		\toprule
    		Method & $p$-hierarchy & $ \Delta t $ & Runtime (s) & $N_\mathrm{ptc}^\mathrm{avg}$ &  $ N_\mathrm{gmres}^\mathrm{avg} $ &   Speedup\\
    		\midrule
    		$p$MG & $ p\{3\text{-}1\text{-}0\} $  & 0.025 & \num{9163} & 96.8 & -- &  1\\
    		$ p $MG &$ p\{3\text{-}2\text{-}1\text{-}0\} $  & 0.025 & \num{9315} & 75 & --  & 0.98\\
    		JFNK-EJ & -- & 0.025 &  \num{4053} & 7.3 & 357.2 &  2.26\\
    		\bottomrule
    	\end{tabular}
    
    \end{table}

	\begin{table}[tbhp]
    	\centering
    	\caption{$ p $MG preconditioner vs. EJ preconditioner for ESDIRK2 for viscous flow over a cylinder at $ \text{Ma}_{\infty} =0.01 $. $ N_\mathrm{ptc}^\mathrm{avg} $ is the averaged number of pseudo time iterations   per RK stage and $ N_\mathrm{gmres}^\mathrm{avg} $ is the averaged number of GMRES iterations per RK stage.}
    	\label{tab:ESDIRK-JFNK-EJ-PMG}
    	\begin{tabular}{cccccccc}
    		\toprule
    		Method & $p$-hierarchy & $\Delta t$ & Runtime (s) & $N_\mathrm{ptc}^\mathrm{avg}$ & $N_\mathrm{gmres}^\mathrm{avg}$ & $K_\mathrm{dim}$ & Speedup\\
    		\midrule
    		$p$MG & $p\{3\text{-}1\text{-}0\}$ & 0.025 & \num{9163} & 96.8 & -- & -- &1\\
    		JFNK-EJ & -- & 0.025 & \num{4053} & 7.3 & \num{357.2} & 30 & 2.26 \\
    		JFNK-EJ & -- & 0.05 & \num{5227} & 8.0 & \num{861.8} & 100 & 1.75\\		
    		JFNK-$p$MG & $p\{3\text{-}1\}$ & 0.05 & \num{2008} & 10.3 & 43.9 & 30 & 4.56\\
    		\bottomrule
    	\end{tabular}
	
	\end{table}

 	\begin{table}[tbhp]
    	\centering
    	\caption{$ p $MG preconditioner vs. EJ preconditioner for ROW2 for viscous flow over a cylinder at $ \text{Ma}_{\infty} =0.01 $. $ K_\mathrm{dim} $ is the maximum dimension of the Krylov subspace. $ N_\mathrm{gmres}^\mathrm{avg} $ is the averaged number of GMRES iterations per RK stage.}
    	\label{tab:ROW-JFNK-EJ-PMG}
    	\begin{tabular}{ccccccc}
    		\toprule
    		Method & $p$-hierarchy & $\Delta t$ & Runtime (s) & $K_\mathrm{dim}$ & $N_\mathrm{gmres}^\mathrm{avg}$ & Speedup\\
    		\midrule
    		$p$MG & $\{3\text{-}1\text{-}0\} $ & 0.025 & \num{9163} & -- & -- & 1 \\
    		GMRES-EJ & -- & 0.05 & \num{3390} & \num{200} & \num{293.7} & 2.70 \\
    		GMRES-$p$MG & $p\{3\text{-}1\}$ & 0.05 & \num{1721} & 30  & 30.8 & 5.32\\
    		GMRES-$p$MG & $p\{3\text{-}0\}$ & 0.05 & \num{5257} & 200 & 102.1 & 1.74 \\
    		GMRES-$p$MG& $ p\{3\text{-}2\text{-}1\} $ & 0.05 & \num{1975} & 30 & 23.6 & 4.64 \\
    		\bottomrule
    	\end{tabular}
	\end{table}
 \clearpage  	
\subsubsection{Unsteady flow over SD7003 airfoil}
\begin{figure}	
		\centering
		\begin{subfigure}[]{0.48\textwidth}
			\centering
			\includegraphics[width=0.99\linewidth]{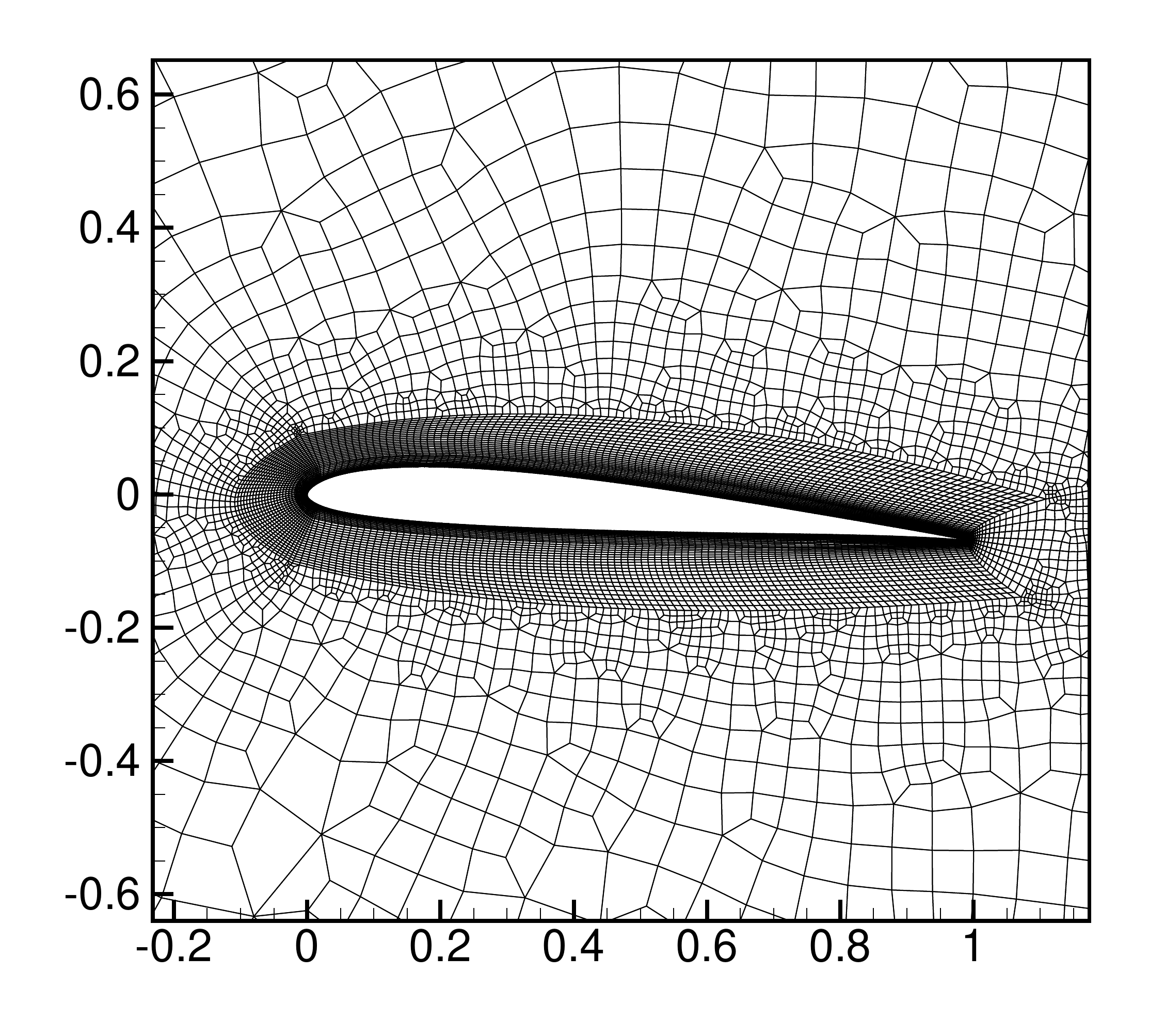}
			\caption{Overview of near-wall mesh}\label{fig:sd7003_nearwall}		
		\end{subfigure}
		\quad
		\begin{subfigure}[]{0.48\textwidth}
			\centering
			\includegraphics[width=0.99\linewidth]{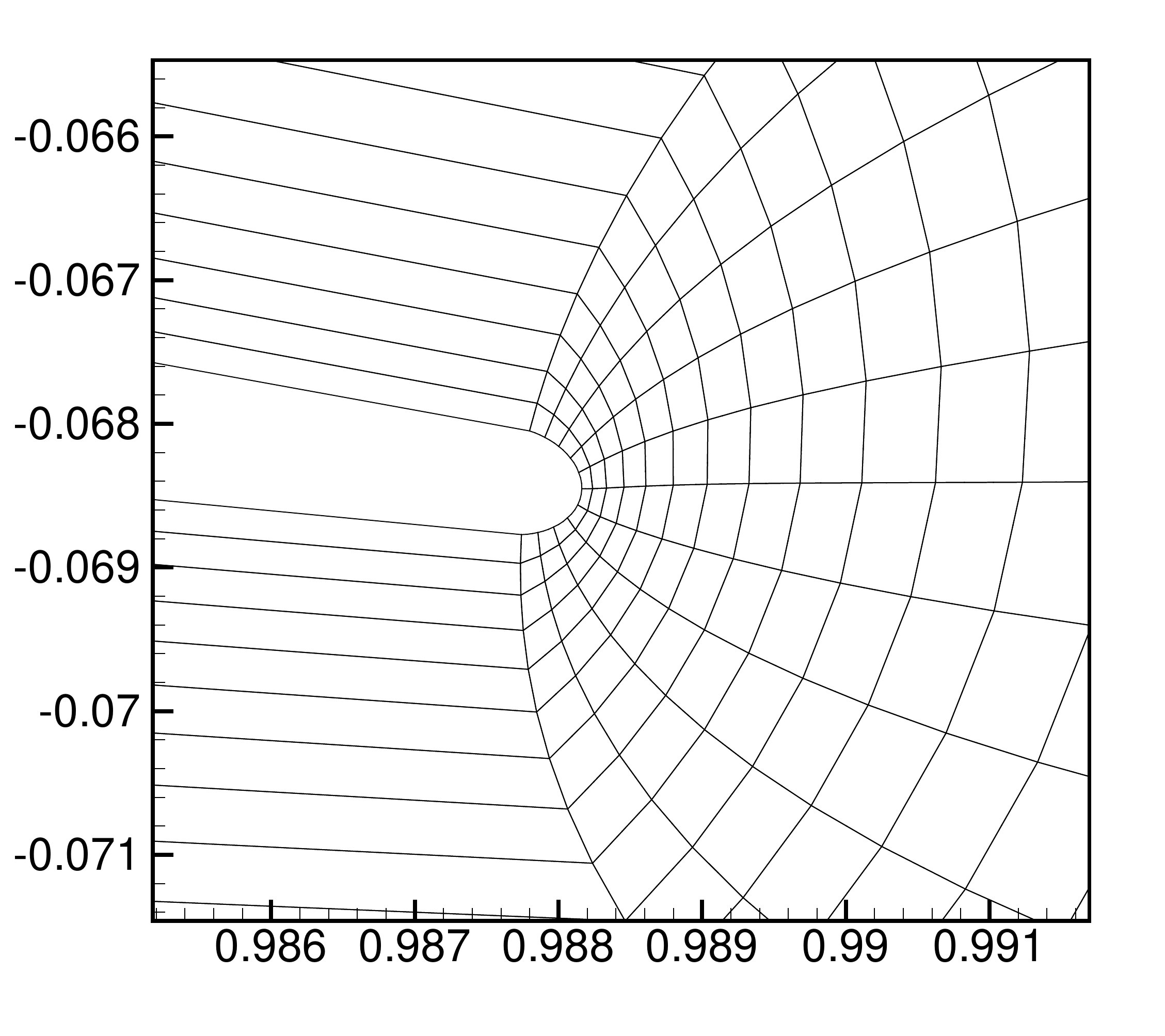}
			\caption{Elements near the trailing edge}\label{fig:sd7003_trailing}
		\end{subfigure}
	\quad
	\begin{subfigure}[]{0.48\textwidth}
		\centering
		\includegraphics[width=0.99\linewidth]{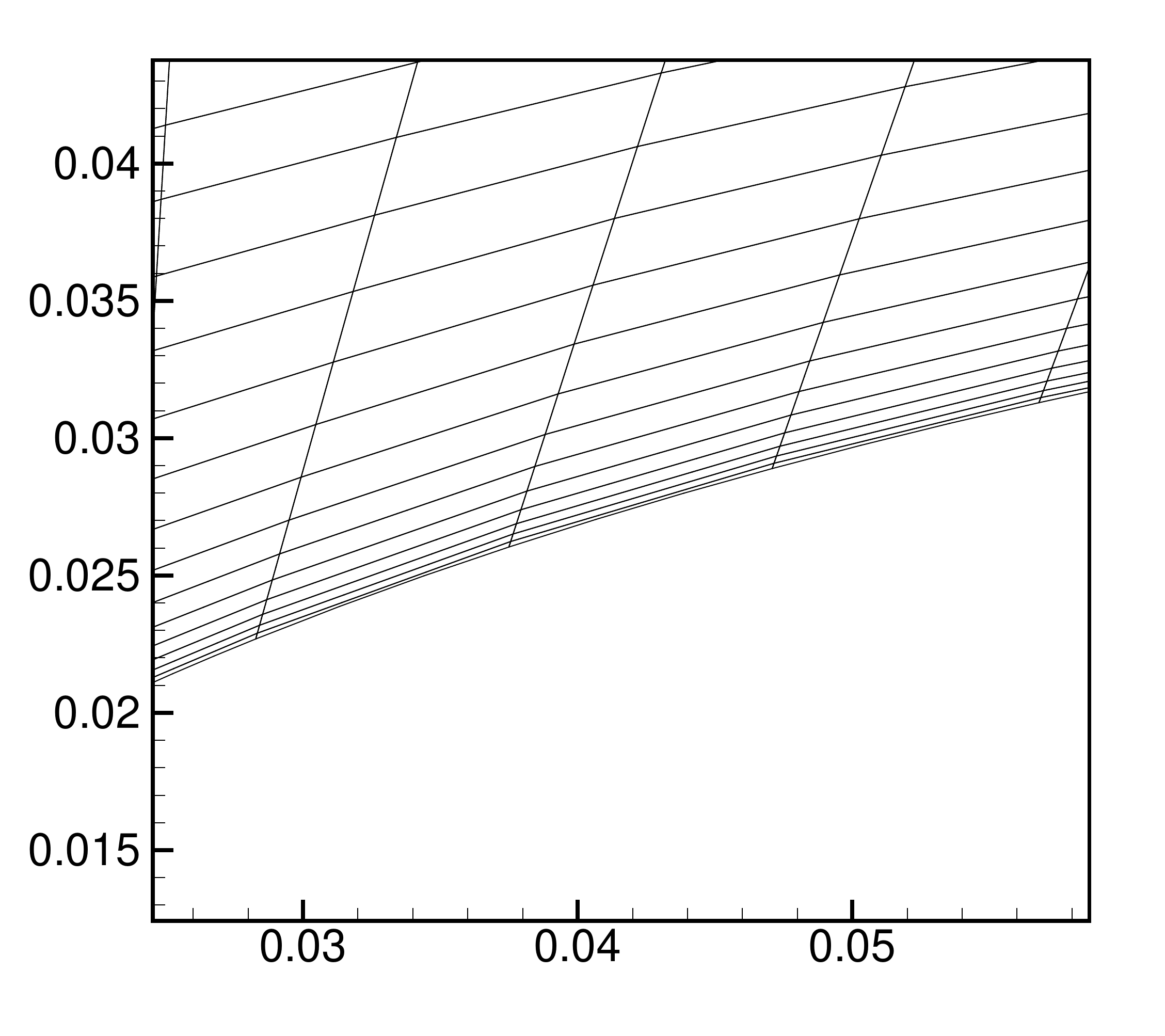}
		\caption{Elements near the leading edge}\label{fig:sd7003_leading}
	\end{subfigure}
		\caption{Near-wall mesh for 2D viscous flow over SD7003 airfoil at $ \text{Re}_c=10^5 $.}
		\label{fig:sd7003_mesh}
	\end{figure}

\begin{figure}[]
    	\centering
    	\adjustbox{width=0.6\linewidth,valign=b}{\includegraphics{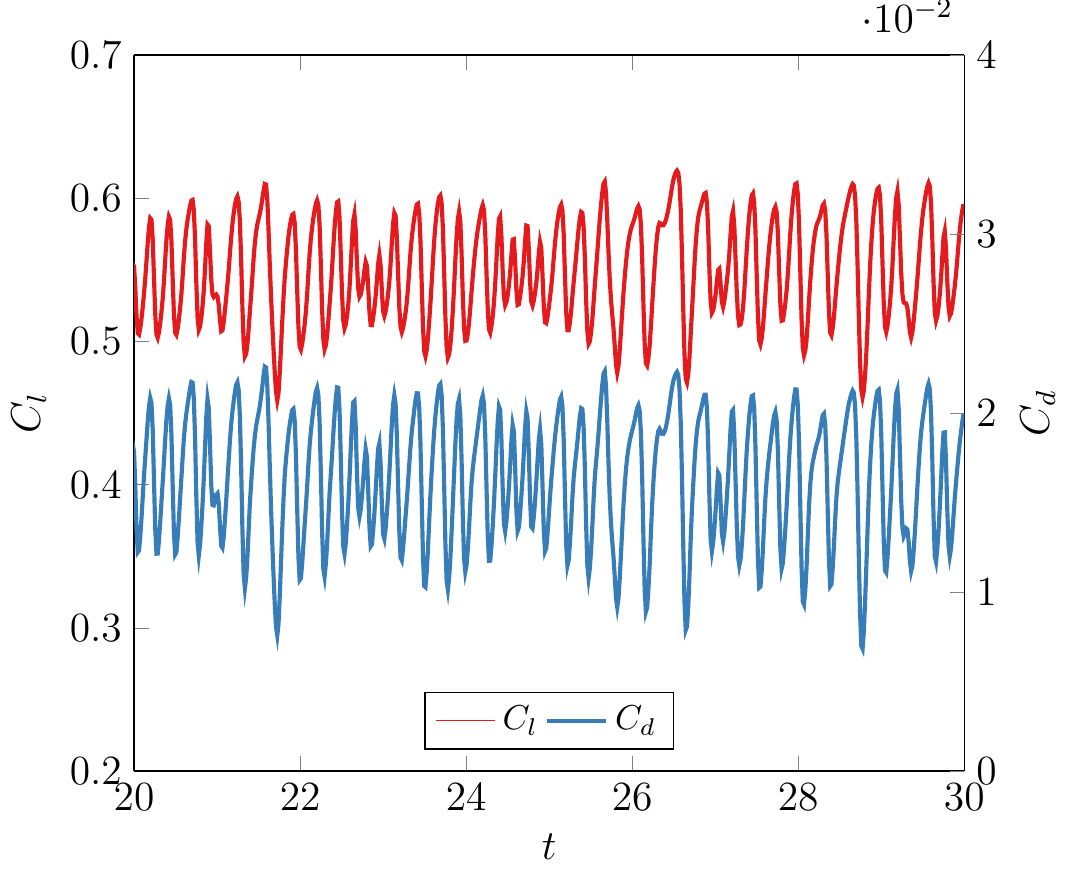}}
    	\caption{\label{fig:sd7003_force}Force histories of 2D vicous flow over SD7003 airfoil.}
    \end{figure}

\begin{figure}	
    	\centering
    	\begin{subfigure}[]{0.48\textwidth}
    		\centering
    		\includegraphics[width=0.99\linewidth]{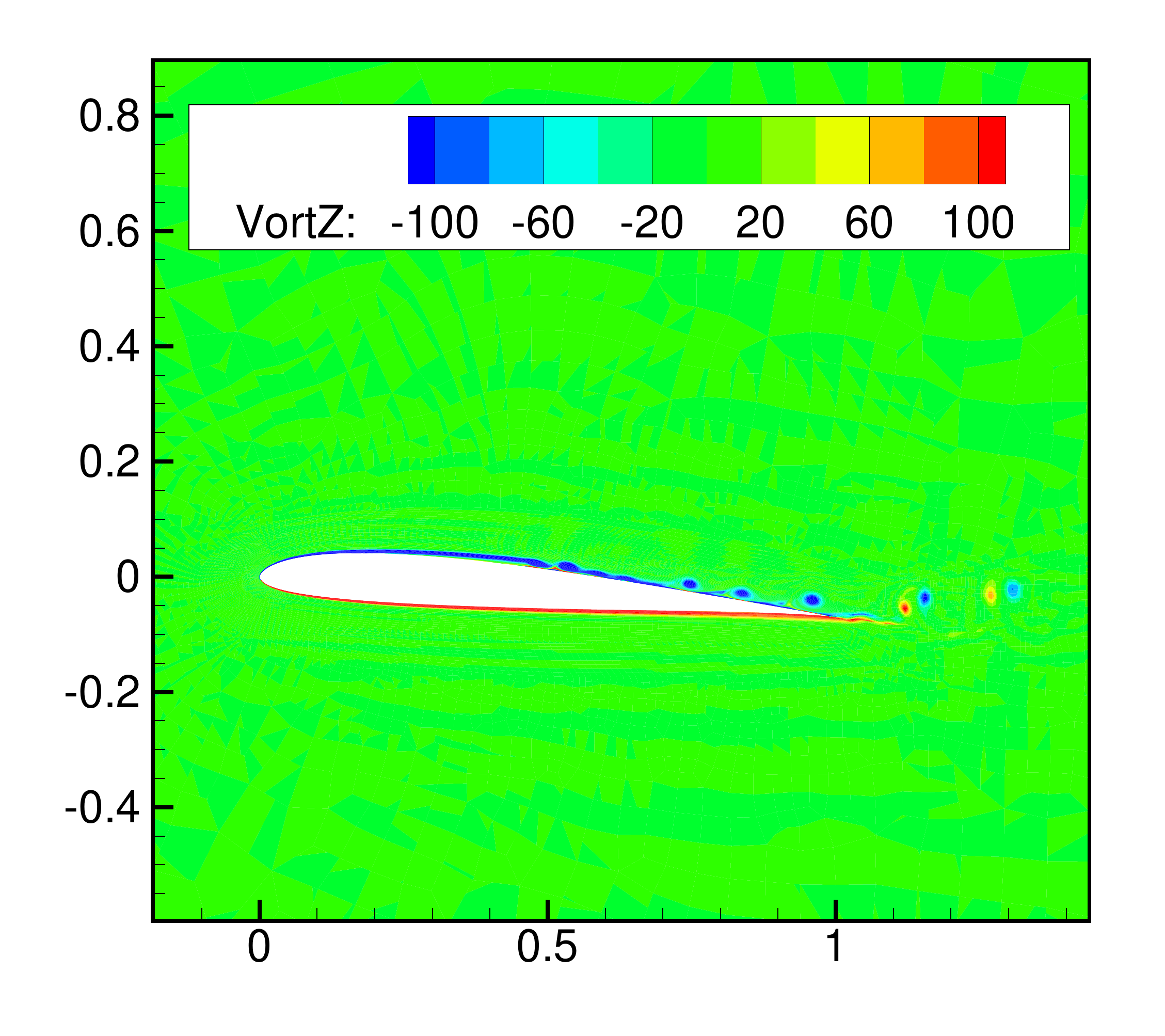}
    		\caption{$z$-Vorticity contour}\label{fig:sd7003_vort}		
    	\end{subfigure}
    	\quad
    	\begin{subfigure}[]{0.48\textwidth}
    		\centering
    		\includegraphics[width=0.99\linewidth]{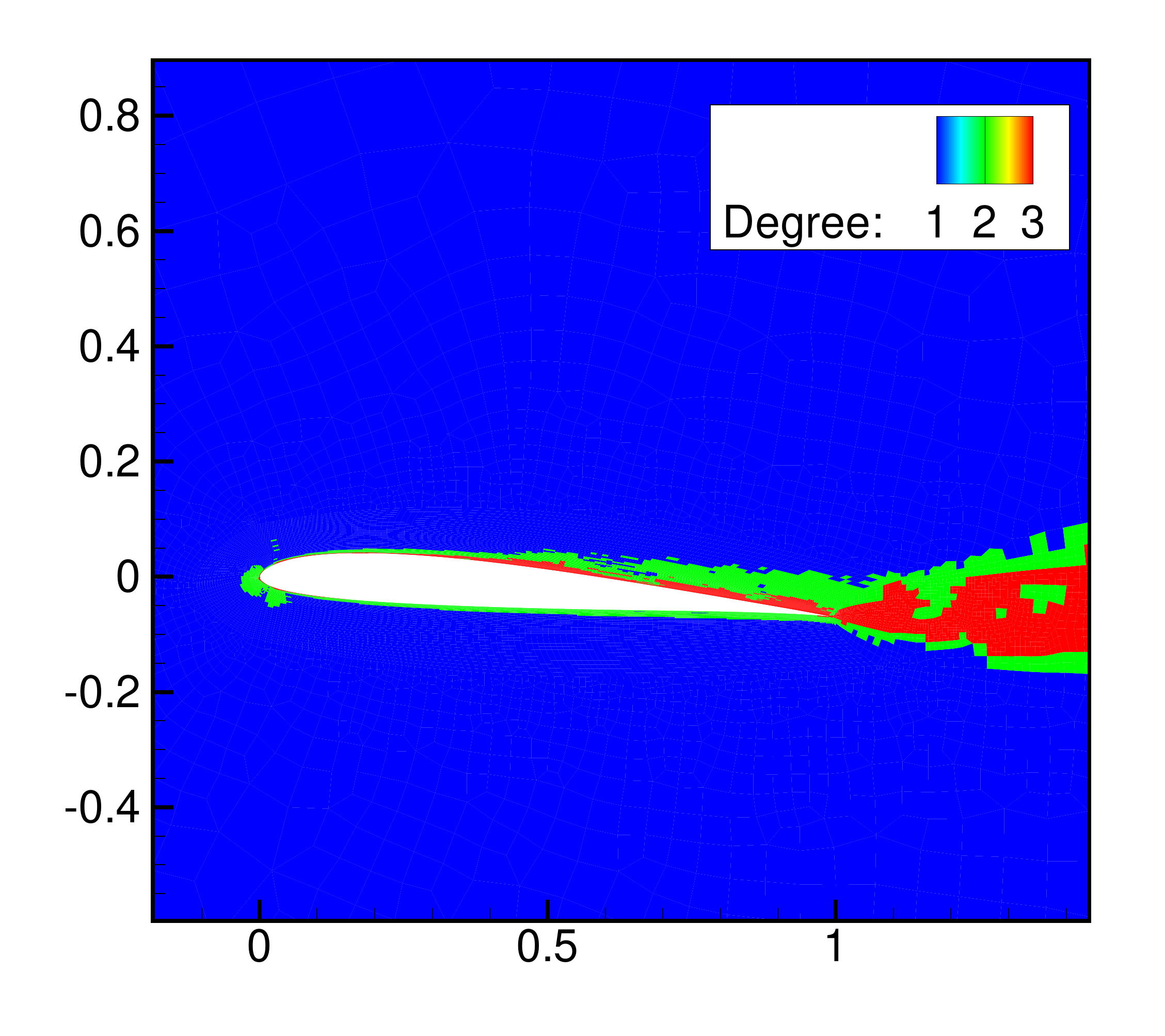}
    		\caption{polynomial disctribution contour}\label{fig:sd7003_order}
    	\end{subfigure}
    	
    	\caption{Instantaneous contours of 2D viscous flow over SD7003 airfoil.}
    	\label{fig:sd7003_contour}
    \end{figure}

	We further examined the performance of $ p $MG preconditioners when dealing with high-aspect-ratio elements clustered to wall boundaries.  For the incoming flow, the Reynolds number based on chord length $ c=1 $ is $ \text{Re}_c = 10^{5} $, Mach number is $ 0.1 $, and the angle of attack is $ 4^\circ $. The nominal height of first layer of elements in vicinity of the wall is $ 2\times 10^{-4}c $ which gives a prior estimation of $ y^+ =1$. In Figure~\ref{fig:sd7003_nearwall}, we illustrate mesh in near wall region; the presence of the small trailing edge, which has a radius of $ 4\times 10^{-4}c $, leads to smaller elements in this region as shown in Figure~\ref{fig:sd7003_trailing}; the mesh generator can only project the linear mesh on to the curved geometry for the first layer of elements to obtain quadratic quadrilateral elements which further worsen the CFL condition in these slim elements as shown in Figure~\ref{fig:sd7003_leading}. The potential of $ p $-adaptation methods on CPU platforms has been demonstrated well for external wall-bounded transitional flows in~\cite{wang2020dynamically}. In this subsection, we employed the $ p $-adaptive $ \mathbb{P}^3 $ FR flux reconstruction for spatial discretization to evaluate $ p $MG methods for non-uniformly distributed polynomial degrees. ESRIRK2 with JFNK-$ p $MG was employed for time integration to run the simulation for $ 30 $ convective time units   to obtain the initial conditions for our further numerical experiments. Force histories are shown in Figure~\ref{fig:sd7003_force}. The averaged  lift and drag coefficients in $t\in(20,30]$ are $C_d=0.0164$ and $C_l=0.5501$. Experimental $C_l$ from~\cite{selig1997summary} is around 0.6. In the transitional regime, a simulation of 2D viscous flow is not sufficient to capture the separation-reattachment phenomenon and current numerical simulation underestimates the lift coefficient by 8.3\%.  Instantaneous contours of vorticity and polynomial degree distribution are demonstrated in Figure~\ref{fig:sd7003_contour}. The $ p $-refined region is consistent with high vorticity region.
	
\begin{figure}[]
    	\centering
    	\adjustbox{width=0.6\linewidth,valign=b}{\includegraphics{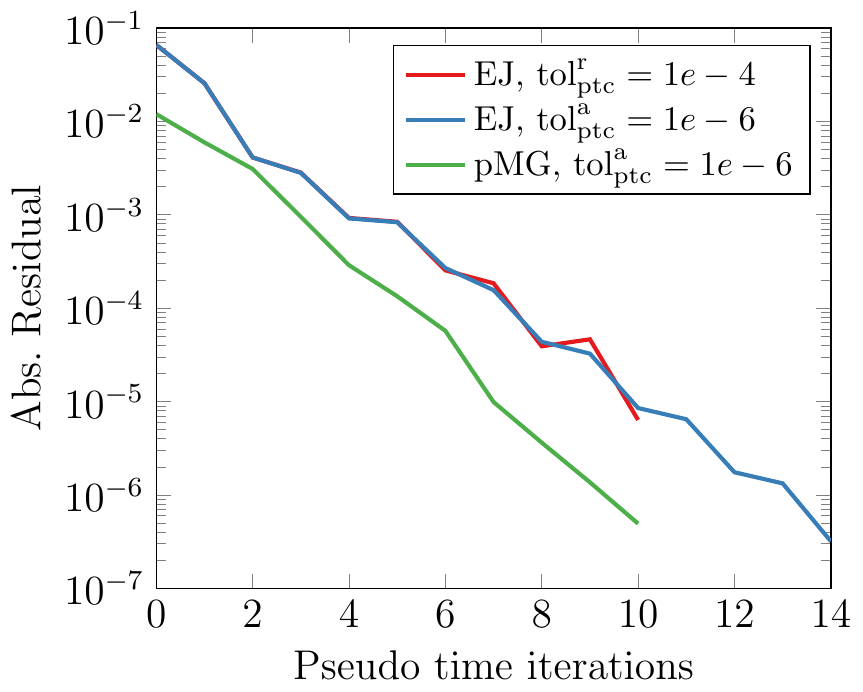}}
    	\caption{\label{fig:sd_res_his}Typical residual history in one stage of ESDIRK4 for 2D viscous flow over SD7003 airfoil.}
    \end{figure}

	We used both EJ and $ p $MG preconditioners for GMRES employed in ESDIRK4 and ROW4 for our comparison. $ \Delta t=0.005 $ was used to resume the simulation for another $ 400 $ time steps. $ \Delta \tau =0.0025 $ was used for pseudo transient continuation in both ESDIRK4 and preconditioner of GMRES.  The $ p $MG preconditioner used in this numerical experiment employed a 2-level polynomial hierarchy and $ n\{2\text{-}2\}  $ for smoothing. In Table~\ref{tab:sd7003_esdirk4}, we present the statistics of  nonlinear/linear solvers used in ESDIRK4. In this table, GMRES-EJ and GMRES-$ p $MG used absolute tolerance $ \mathrm{tol}_\mathrm{ptc}^\mathrm{a} = 10^{-6} $ as the convergence tolerance for PTC and GMRES-EJ$ ^* $ used the relative tolerance $ \mathrm{tol}_\mathrm{ptc}^\mathrm{r} =10^{-4}$. As shown in Figure~\ref{fig:sd_res_his},  $ p $MG preconditioner made the initial residual much smaller than EJ preconditioner did after the first pseudo time iteration. Therefore, using $ \mathrm{tol}_\mathrm{ptc}^\mathrm{a} = 10^{-6} $ as the tolerance is important for a fair comparison. In this case, $ p $MG preconditioner can achieve speedup of 1.20 compared to EJ preconditioner with remarkably smaller Krylov subspace dimension and fewer pseudo time iterations. In Table~\ref{tab:sd7003_row4}, we present the statistics of ROW4 using different preconditioners. $ \mathrm{tol}_\mathrm{gmres}^\mathrm{r}=10^{-6} $ was used to preserve the nominal order of accuracy. We have observed that a speedup of 1.49 via employing the $ p $MG preconditioner with $ p\{3\text{-}2\} $ instead of EJ. Note that insufficient smoothing on intermediate sublevel when using $ p\{3\text{-}1\} $ as the $ p $-hierarchy worsened the performance of $ p $MG. 
	
	For this type of problem, with a $ p $MG preconditioner, ROW4 is more efficient than ESDIRK4 when the nominal order of accuracy is to be preserved via driving the residual sufficiently small. More importantly, ROW4 tends to benefit more from employing such preconditioner in terms of both speedup and reduction of the absolute number of dimension of Krylov subspace.
	
	We are aware that in LES, sufficient resolution in all directions is warranted such that the aspect ratio is less likely to be as huge as it would be in RANS.
	However, for practical applications with more complicated geometries, it is common to observe highly distorted elements or elements which are significantly smaller than what $y^+$ requires at locations where the curvature is large, such as the trailing edge of the SD7003 airfoil studied here. The $p$MG preconditioner developed here has the significance of improving the robustness of the CFD solver for bad quality meshes.

	\begin{table}
		\centering
		\caption{GMRES with $ p $MG preconditioner vs.~EJ preconditioner for ESDIRK4 for viscous flow over SD7003 airfoil with high-aspect-ratio elements in vicinity of wall. Note that the $ p $-hierarchy is only the one on elements which have $ \mathbb{P}^3 $ polynomials.}
		\label{tab:sd7003_esdirk4}
		\begin{tabular}{cccccccc}
			\toprule
			Precond. & $p$-hierarchy & $\Delta t$ & Runtime (s) & $ N_\mathrm{ptc}^\mathrm{avg} $ & $N_\mathrm{gmres}^\mathrm{avg}$ & $K_\mathrm{dim}$ & Speedup\\
			\midrule
			EJ & -- & 0.005 & \num{23363} & 42.8 & \num{865.8} & 30 & 1 \\
			EJ$^*$ & -- & 0.005 & \num{17307} & 32.3 & \num{632.4} & 30 & 1.35 \\
			$p$MG & $p\{3\text{-}2\}$ & 0.005 & \num{19396} & 26.3 & 79.1 & 5 & 1.20 \\
			\bottomrule
		\end{tabular}	
	\end{table}

	\begin{table}
		\centering
		\caption{GMRES with $p$MG preconditioner vs.~EJ preconditioner for ROW4 for viscous flow over SD7003 airfoil with high-aspect-ratio elements in vicinity of wall. Note that the $ p $-hierarchy is only the one on elements which have $ \mathbb{P}^3 $ polynomials.}
		\label{tab:sd7003_row4}
		\begin{tabular}{ccccccc}
			\toprule
			Precond. & $p$-hierarchy & $\Delta t$ & Runtime (s) & $K_\mathrm{dim}$ & $N_\mathrm{gmres}^\mathrm{avg}$ & Speedup\\
			\midrule
			EJ & -- & \num{0.005} & \num{18463} & \num{150} & \num{185.4} &1 \\		
			$p$MG & $ p\{3\text{-}2\} $& \num{0.005} & \num{12392} & 30 & 22.5 & 1.49 \\
			$p$MG & $ p\{3\text{-}1\} $& \num{0.005} & \num{34876} & 30 & 45.4 & 0.53 \\
			\bottomrule
		\end{tabular}
		
	\end{table}

	
\section{Conclusions and future work}\label{sec:conclusion}
	
	We have investigated the performance of $ p $MG when $ p $-adaptation is coupled. 
	We start our discussion of $ p $MG by analyzing its performance as a nonlinear solver for steady problems. It is demonstrated that employment of a stronger smoother at coarsest $ p $-sublevel is critically important for efficiency. Moreover, insufficient smoothing on intermediate $ p $-sublevel would even introduce instability into $ p $MG nonlinear solvers. We  hypothesise that in PTC, as the residual drops, the frequency of error modes that dictates the convergence rate get higher and higher, which have been demonstrated through numerical experiment. This is particularly important for LES since the initial guess in PTC used in implicit time integration is not distance from the converged solution.
	
	We propose to use pseudo transient continuation in the preconditioning procedure to construct nonlinear $ p $MG preconditioners for matrix-free GMRES used in implicit time integration of unsteady problems. 
	We have demonstrated that for stiff systems, JFNK with even just an element-Jacobi preconditioner can significantly outperform $ p $MG nonlinear solver by a factor over 2. JFNK with $ p $MG preconditioner can be around 2 times faster for low-Mach-number flows and 1.5 times faster for anisotropic meshes than JFNK-EJ. More importantly,  similar to what has been observed in simulation of steady problems, insufficient smoothing on intermediate $ p $-sublevel can greatly worsen the performance of $ p $MG which has been overlooked in literature. We genuinely recommend a $ \{p_0\text{-}p_0/2\} $ polynomial hierarchy for $ p $MG preconditioners when used for unsteady problems, or even $ \{p_0\text{-}(p_0-1)\} $ when RAM usage is affordable to use a MBNK smoother at the bottom. Additionally, it is observed that with a decent  preconditioner, namely $ p $MG in this study, ROW is consistently more efficient than ESDIRK when the residual is required to drop to small enough to preserve the order of accuracy. In particular, the reduction of the absolute value of the dimension of Krylov subspace for ROW is more significant than that for ESDIRK. 
	
	Researchers recently started to look into employing JFNK methods on GPUs for LES~\cite{jourdan2021matrix}. We are aware that storing the element-Jacobi smoother of the finest $ p $-sublevel is very overwhelming for such hardware, especially when the polynomial degree is high. For JFNK-$ p $MG on GPU, we would like to explore the feasibility of using a first order exponential  time integrator~\cite{li2018exponential} as the smoother instead of the EJ smoother and the matrix-based smoother. The overall RAM usage can potentially be minimised without the presence of an element-Jacobi smoother on the finest level when the Krylov subspace dimension is limited and such smoother would be able to march in pseudo time with large strides since the stability can be improved with a $p$MG configuration. Investigation on such topic will be our future work.
	
\section*{Acknowledgement}
The first author would like thank Tarik Dzanic for proofreading this manuscript.

\clearpage
\bibliography{ptc-preconditioner}

\section*{Statements \& Declarations}
\subsection*{Funding}
The authors declare that no funds, grants, or other support were received during the preparation of this manuscript.

\subsection*{Competing Interests}
The authors have no relevant financial or non-financial interests to disclose.

\subsection*{Data Availability}
The datasets generated during and/or analysed during the current study are available from the corresponding author on reasonable request.
\end{document}